\documentclass[11pt,reqno]{amsart}

\usepackage{amsthm}
\usepackage{amssymb}
\usepackage{latexsym}
\usepackage{multicol}
\usepackage{verbatim,enumerate}
\usepackage{hyperref}
\usepackage{nicematrix}
\usepackage{amsmath, amscd}
\usepackage{mathrsfs}
\usepackage[all,cmtip]{xy}
\usepackage{soul}

\advance\textwidth by 1.2in \advance\oddsidemargin by -.6in \advance\evensidemargin by -.6in
\usepackage{tikz}
\usetikzlibrary{decorations.markings}
\tikzstyle{vertex}=[ circle, fill, draw, inner sep=0pt, minimum size=4pt,]
\tikzstyle{edge}= [thick]

\newtheorem*{cor}{Corollary}%[section]
\newtheorem*{lem}{Lemma}
\newtheorem*{prop}{Proposition}

\theoremstyle{definition} 
\theoremstyle{definition}
\newtheorem{thm}{Theorem}
\newtheorem*{thm*}{Theorem}

\newtheorem*{rem}{Remark}

\newenvironment{pf}{\proof}{\endproof}
\newcounter{cnt}
\newenvironment{enumerit}{\begin{list}{{\hfill\rm(\roman{cnt})\hfill}}{%
\settowidth{\labelwidth}{{\rm(iv)}}\leftmargin=\labelwidth%
\advance\leftmargin by \labelsep\rightmargin=0pt\usecounter{cnt}}}{\end{list}} \makeatletter
\def\mydggeometry{\makeatletter\dg@YGRID=1\dg@XGRID=20\unitlength=0.003pt\makeatother}
\makeatother \theoremstyle{remark}

\numberwithin{equation}{section}

\newcommand{\nc}{\newcommand}
\newcommand{\rnc}{\renewcommand}

\nc{\cal}{\mathcal} \nc{\goth}{\mathfrak} \rnc{\bold}{\mathbf}

\renewcommand{\Bbb}{\mathbb}
\nc\bomega{{\mbox{\boldmath $\omega$}}} \nc\bpsi{{\mbox{\boldmath $\Psi$}}}
 \nc\sing{{\rm sing}}
 \nc\balpha{{\mbox{\boldmath $\alpha$}}}
 \nc\bbeta{{\mbox{\boldmath $\beta$}}}
 \nc\bpi{{\mbox{\boldmath $\pi$}}}
  \nc\bpis{{\mbox{\boldmath \scriptsize$\pi$}}}
 \nc\bullets{{\mbox{\scriptsize $\bullet$}}}
 \nc\bvarpis{{\mbox{\boldmath \scriptsize$\varpi$}}}
  \nc\bvarpi{{\mbox{\boldmath $\varpi$}}}

\nc\bepsilon{{\mbox{\boldmath $\epsilon$}}}

  \nc\bomegas{{\mbox{\boldmath\scriptsize $\omega$}}}
  \nc\bepsilons{{\mbox{\boldmath \scriptsize$\epsilon$}}}
\nc{\spi}{{\rm sp}}
\nc\hlien{\hat{\lie n}^+}
  \nc\btaus{{\mbox{\boldmath \scriptsize$\tau$}}}\nc\bxi{{\mbox{\boldmath $\xi$}}}
\nc\bmu{{\mbox{\boldmath $\mu$}}} \nc\bcN{{\mbox{\boldmath $\cal{N}$}}} \nc\bcm{{\mbox{\boldmath $\cal{M}$}}} \nc\blambda{{\mbox{\boldmath
$\lambda$}}}%\nc\mathbb Nu{{\mbox{\boldmath $\nu$}}}
\nc\btau{{\mbox{\boldmath
$\tau$}}}

\newcommand{\lie}[1]{\mathfrak{#1}}

\makeatletter
\def\section{\def\@secnumfont{\mdseries}\@startsection{section}{1}%
  \z@{.7\linespacing\@plus\linespacing}{.5\linespacing}%
  {\normalfont\scshape\centering}}
\def\subsection{\def\@secnumfont{\bfseries}\@startsection{subsection}{2}%
  {\parindent}{.5\linespacing\@plus.7\linespacing}{-.5em}%
  {\normalfont\bfseries}}
\makeatother

 \nc{\Hom}{\operatorname{Hom}}
  \nc{\mode}{\operatorname{mod}}
\nc{\End}{\operatorname{End}} \nc{\wh}[1]{\widehat{#1}} \nc{\Ext}{\operatorname{Ext}}
 \nc{\ch}{\operatorname{ch}} \nc{\ev}{\operatorname{ev}}
\nc{\Ob}{\operatorname{Ob}} \nc{\soc}{\operatorname{soc}} \nc{\rad}{\operatorname{rad}} \nc{\head}{\operatorname{head}}

\def\mult{\operatorname{mult}}

 \nc{\Cal}{\cal} \nc{\Xp}[1]{X^+(#1)} \nc{\Xm}[1]{X^-(#1)}
\nc{\on}{\operatorname} \nc{\Z}{{\bold Z}} \nc{\J}{{\cal J}} \nc{\C}{{\bold C}} \nc{\Q}{{\bold Q}}

\nc{\N}{{\Bbb N}} \nc\boa{\bold a} \nc\bob{\bold b} \nc\boc{\bold c} \nc\bod{\bold d} \nc\boe{\bold e} \nc\bof{\bold f} \nc\bog{\bold g}
\nc\boh{\bold h} \nc\boi{\bold i} \nc\boj{\bold j} \nc\bok{\bold k} \nc\bol{\bold l} \nc\bom{\bold m} \nc\bon{\bold n} \nc\boo{\bold o}
\nc\bop{\bold p} \nc\boq{\bold q} \nc\bor{\bold r} \nc\bos{\bold s} \nc\boT{\bold t} \nc\boF{\bold F} \nc\bou{\bold u} \nc\bov{\bold v}
\nc\bow{\bold w} \nc\boz{\bold z} \nc\boy{\bold y} \nc\ba{\bold A} \nc\bb{\bold B} \nc\bc{\mathbb C} \nc\bd{\bold D} \nc\be{\bold E} \nc\bg{\bold
G} \nc\bh{\bold H} \nc\bi{\bold I} \nc\bj{\bold J} \nc\bk{\bold K} \nc\bl{\bold L} \nc\bm{\bold M}  \nc\bo{\bold O} \nc\bp{\bold
P} \nc\bq{\bold Q} \nc\br{\bold R} \nc\bs{\bold S} \nc\bt{\bold T} \nc\bu{\bold U} \nc\bv{\bold V} \nc\bw{\bold W} \nc\bx{\bold
x} \nc\KR{\bold{KR}} \nc\rk{\bold{rk}} \nc\het{\text{ht }}
\nc\bz{\mathbb Z}
\nc\bn{\mathbb N}
\nc\us{\underline \bos}
\nc\uS{ \bs_{{\rm alt}}}
\nc\pr{\rm pr}

\nc\toa{\tilde a} \nc\tob{\tilde b} \nc\toc{\tilde c} \nc\tod{\tilde d} \nc\toe{\tilde e} \nc\tof{\tilde f} \nc\tog{\tilde g} \nc\toh{\tilde h}
\nc\toi{\tilde i} \nc\toj{\tilde j} \nc\tok{\tilde k} \nc\tol{\tilde l} \nc\tom{\tilde m} \nc\ton{\tilde n} \nc\too{\tilde o} \nc\toq{\tilde q}
\nc\tor{\tilde r} \nc\tos{\tilde s} \nc\toT{\tilde t} \nc\tou{\tilde u} \nc\tov{\tilde v} \nc\tow{\tilde w} \nc\toz{\tilde z} \nc\woi{w_{\omega_i}}
\nc\chara{\operatorname{Char}}
\def\op{\operatorname{op}}
\nc\Sym{{\rm Sym}}
\nc\In{\mathbb I_n}
\nc\Fn{\mathscr F_n}
\begin{document}

\title[]{The Zelevinsky involution, co-standard modules for quantum affine algebras and interval symmetric functions}
\author{Matheus Brito}
\address{Departamento de Matematica, UFPR, Curitiba - PR - Brazil, 81530-015}
\email{mbrito@ufpr.br}
\thanks{M.B. was partially supported by CNPq grants 405793/2023-5 and 300178/2026-3}
\author{Vyjayanthi Chari}
\address{Department of Mathematics, University of California, Riverside, 900 University Ave., Riverside, CA 92521, USA}
\email{chari@math.ucr.edu}
\thanks{V.C. was partially supported by a travel grant from the Simons Foundation.}

\begin{abstract}
Let $\mathscr F_n$ be the subcategory of finite-dimensional representations of the quantum loop algebra of $\lie{sl}_{n+1}$ whose
simple modules are indexed by multisegments consisting of intervals $[a,b]$, where $a,b$ are integers and $0<b-a\le n$. We define
the co-standard module associated to a multisegment $\bos$ to be the tensor product of the Kirillov--Reshetikhin modules,
associated to the first node of the Dynkin diagram of $\lie{sl}_{n+1}$, given by the segments of $\bos$. Under a suitable (and
necessary) restriction on $n$, we show that its socle is isomorphic to the simple module given by the Zelevinsky dual $\bos^\dagger$
of $\bos$. We then define an involution $\Phi$ of the Grothendieck ring $\cal K(\mathscr F_n)$ and prove that it maps the class of a simple
module to the class of the simple module given by its Zelevinsky dual. We also show that on suitable subcategories $\Phi$ is
induced by an exact autoequivalence reversing tensor products, and we give a $t$-analogue of the multiplicity identity in the
quantum Grothendieck ring.
\\\\
Motivated by the determinantal formula for alternating snake modules, we realize $\cal K(\mathscr F_n)$ explicitly inside a
polynomial ring in variables indexed by intervals, and, in the spirit of Macdonald's ninth variation of Schur functions, we
show that the fundamental modules correspond to analogues of the elementary symmetric functions and the
Kirillov--Reshetikhin modules to analogues of the complete homogeneous symmetric functions. This is compatible with $\Phi$ and
the classical involution $\omega$, and we prove that alternating snake modules give a family of Schur positive symmetric
functions extending the skew Schur functions.
\end{abstract}

\maketitle

\section*{Introduction}

For $p$-adic general linear groups, Zelevinsky's classification of irreducible representations by multisegments comes with
two families of standard modules, products of the segment representations $Z(\Delta)$ and of the representations $L(\Delta)$,
and with an involution exchanging them which maps irreducible representations to irreducible representations; on
multisegments this involution is computed by the algorithm of M\oe glin and Waldspurger \cite{MW86}. In this paper we
develop the analogous theory for finite-dimensional representations of the quantum loop algebra $\widehat\bu_n$ of $\lie{sl}_{n+1}$,
using the methods based on $R$-matrices introduced by Kashiwara, Kim, Oh and Park \cite{KKOP}.
\\\\
We work in the tensor subcategory $\mathscr F_n$ of Hernandez and Leclerc \cite{HL10}, whose simple objects are indexed by the
same multisegments, now of intervals $[a,b]$ with $0<b-a\le n$, and whose Grothendieck ring is the polynomial ring in the
classes of the fundamental modules \cite{FR99}. Several families of modules in $\mathscr F_n$ have counterparts on the $p$-adic
side, which we use as a guide; for instance the snake modules of Mukhin and Young \cite{MY12} correspond to the ladder
representations of Lapid and M\'inguez \cite{LM14}.
\\\\
In $\mathscr F_n$ the standard module $W(\bomega_\bos)$ attached to a multisegment $\bos$ is a tensor product of fundamental modules,
one for each segment of $\bos$. We define the co-standard module $M(\bos)$ to be the tensor product of the
Kirillov--Reshetikhin modules associated to the first node and to the segments of $\bos$, and we consider the ring involution
$\Phi$ of the Grothendieck ring sending the class of a fundamental module to the class of the corresponding
Kirillov--Reshetikhin module. Thus $\Phi$ maps standard classes to co-standard classes, and it is the analogue of the involution
$\omega$ of the ring of symmetric functions. Our main result is that co-standard modules have simple head, and that if $n$ is
at least the length $n(\bos)$ of the support of $\bos$ then $M(\bos)$ has simple socle, given by the Zelevinsky dual $\bos^\dagger$ of
$\bos$, and $\Phi$ maps the class of the simple module given by $\bos$ to that given by $\bos^\dagger$ (Theorem \ref{zdualsoc}). The
restriction on $n$ is necessary. As a consequence $\Phi$ exchanges the heads and socles of standard and co-standard modules,
and the multiplicities of simple modules in co-standard modules are those in standard modules, with the labels replaced by
their Zelevinsky duals. We also obtain a $t$-analogue of the identity between
multiplicities in standard and co-standard modules.
\\\\
The theorem also gives information about real and prime modules. Since $\Phi$ preserves irreducibility of tensor products in
our range, a simple module is real (resp. prime) if and only if the module given by the Zelevinsky dual is. This is of
interest because imaginary modules, those which do not correspond to cluster monomials in the categorification of
\cite{HL13a,KKOP}, are hard to come by: the families constructed in \cite{BC23,BM25} all have an even number of segments (but, see
also the recent work \cite{LP26,Sak25}). The Zelevinsky duals of the modules in \cite{BC23,BM25} are usually new imaginary modules,
often with an odd number of segments, and they are prime whenever the original modules are. Moreover there are examples where
neither the multisegment nor its dual is regular, so that these modules are not covered by the criterion of Lapid and
M\'inguez \cite{LM18}.
\\\\
The socle statement is proved directly within $\mathscr F_n$, using the invariants of \cite{KKOP} and their theory of normal
sequences: we realize $M(\bos)$, through a chain of inclusions, inside a tensor product of a normal sequence of real simple
modules, chosen so as to reflect the recursion of the M\oe glin--Waldspurger algorithm. The statement on $\Phi$ then follows
once one knows that $\Phi$ maps simple classes to simple classes, for which we use the generalized Schur--Weyl duality functor of
Kang, Kashiwara and Kim \cite{KKK15}. When $|\bos|$ is small compared with $n$, part of these results could also be transported
from affine Hecke algebras, and hence from $p$-adic groups, through the quantum affine Schur--Weyl duality of \cite{CP96}, which is an equivalence in that range and matches highest $\ell$-weights with multisegments; our theorem
holds whenever $n\ge n(\bos)$, independently of $|\bos|$, where this duality is no longer an equivalence. On the subcategories
$\mathscr F_n(\bou_c)$ one could instead work on the quiver Hecke side, but this would require combining several deep results:
Fujita's theorem that the functor of \cite{KKK15} is an equivalence and coincides in this case with the Schur--Weyl duality of
\cite{CP96} \cite{Fuj}, the identification of quiver Hecke algebras with (completions of) affine Hecke algebras
\cite{BK09,Rou08}, the matching of labellings of simple modules in all degrees, and the description of Kashiwara's involution in
terms of the Zelevinsky involution \cite{LTV99,JL10}. Our arguments use none of these. However, using Fujita's equivalence we show that on these subcategories $\Phi$
is induced by an exact autoequivalence reversing tensor products.
\\\\
Our second theme is a symmetric function interpretation, motivated by the determinantal formula of \cite{BC25a} for
alternating snake modules. These are a class of simple modules which include the snake modules of \cite{MY12} and the modules of the
category $\cal C_1$ of \cite{HL10}. It was shown that the class of a stable alternating snake module is the determinant of an explicit matrix $A(\bos)$
whose entries are classes of fundamental modules or zero. \\\\
The main new ingredient here is an explicit realization of the
Grothendieck ring inside the polynomial ring $\mathbb Z[\mathbb I_n^+]$ in variables $x_{a,b}$ indexed by intervals. The classes of the
fundamental modules and of the Kirillov--Reshetikhin modules become explicit sums of monomials in $x_{a,b}$. Denoting these sums by $e_{a,b}$ and $h_{a,b}$ respectively, we show that they satisfy the identity relating elementary and complete
homogeneous symmetric functions. The ring ${\rm Sym}(\mathbb I_n)$ generated by the $e_{a,b}$ is in the spirit of Macdonald's ninth variation of
Schur functions \cite{Mac92}.\\\\
We relate ${\rm Sym}(\mathbb I_n)$ to the usual ring of symmetric functions, and interpret the main theorem as the analogue
of $\omega(s_\lambda)=s_{\lambda'}$. For snake modules the analogy is literal: the Zelevinsky dual corresponds to transposition of skew
Young diagrams, as observed for ladders in \cite{LM14}.
For regular multisegments whose associated permutation avoids the patterns $3412$ and $4231$, Bittmann \cite{Bit23}
obtained, from the alternating sum formula of Lapid and M\'inguez \cite{LM18}, a determinantal formula for $p$-adic general
linear groups whose matrix also has zeros, determined by a Bruhat interval; it transfers to $\mathscr F_n$ through the
Schur--Weyl duality of \cite{CP96}, and hence only when $|\bos|$ is small compared with $n$. The formula of \cite{BC25a} holds
whenever $n\ge n(\bos)$, and it also applies to alternating snakes outside her class: for instance the prime stable
alternating snake $([0,4],[3,5],[1,3],[0,2])$ is not regular, and $([0,3],[3,7],[6,10],[5,6])$ is regular but its associated
permutation is $4231$.\\\\
For the reader's convenience, we collect the definitions and state the main results in Section 1; the symmetric function
interpretation is developed in Section \ref{lambdasec}, the proofs of the main results are given in Sections
\ref{head}--\ref{phi}, and the autoequivalence and $t$-analogues are discussed in Section \ref{sec:sigma}.

\section{Main Results}\label{sec:main}

We introduce the main objects of study and state the main results of the paper.

\subsection{} As usual $\mathbb C$ (resp. $\mathbb C^\times$, $\mathbb Z$, $\mathbb N$) denotes the set of complex numbers (resp. non-zero complex numbers, integers, positive integers). Given any $r\in\mathbb N$ we let $\Sigma_r$ be the group of permutations of $r$ letters.\\\\
Assume throughout that $q\in \mathbb C^\times$ is not a root of unity.

\subsection{The algebra $\widehat\bu_n$ and the category $\mathscr F_n$} \label{ms}For $n\in \mathbb N$, let $\widehat\bu_n$ be the quantum loop algebra associated to $\lie{sl}_{n+1}$; it is a Hopf algebra and we refer the reader to \cite{CP91,CP94} for precise definitions. \\\\
We shall be interested in a particular subcategory $\mathscr F_n$ of finite-dimensional representations of this algebra.
 Following \cite{BC25b}, let $\mathbb I_n$ be the set of intervals $[a,b]$ with $a,b\in\mathbb Z$ and $0\le b-a\le n+1$, and let $\mathbb I_n^+$ be the subset of intervals $[a,b]$ with $a<b$. We say that a pair of intervals $([a_1, b_1],[a_2,b_2])$ is overlapping  if for some $\epsilon\in\{0,1\}$ we have
 $$a_{1+\epsilon}<a_{2-\epsilon}\le b_{1+\epsilon}<b_{2-\epsilon}.$$ An overlapping pair $([a_1, b_1],[a_2,b_2])$ is connected if in addition we have  $[a_{1+\epsilon}, b_{2-\epsilon}]\in\mathbb I_n$.
 Denote by ${\rm mult}(\mathbb I_n)$ the set of multisegments $\bos=\{[c_1,d_1],\cdots, [c_r,d_r]\}$ with $r\ge 1$ and $[c_s, d_s]\in\mathbb I_n^+$ for $1\le s\le r$. Let $\mu_\bos:\mathbb I_n^+\to\mathbb Z_+$ be the multiplicity function of $\bos$; namely, $\mu_\bos([c,d])$ is the number of times $[c,d]$ appears in $\bos$. Given $\bos_1,\bos_2\in\mult\In$ let $\bos_1+\bos_2$ be the multisegment defined by $\mu_{\bos_1}+\mu_{\bos_2}$.\\\\
 {\em In what follows (unless otherwise stated) we shall assume that a multisegment $\{[c_1,d_1],\cdots, [c_r, d_r]\}$ satisfies $0<d_j-c_j<n+1$ for all $1\le j\le r$. We shall also assume
that it is ordered so that $d_1\ge\cdots\ge d_r$, and if $d_j=d_{j+1}$ then $c_j\ge c_{j+1}$. Set $$n(\bos)= \max\{d_j:1\le j\le r\}-\min\{c_j: 1\le j\le r\},\qquad |\bos|=\sum_{j=1}^r(d_j-c_j).$$}
  Let $\cal I_n$ (resp. $\cal I_n^+$) be the free abelian group (resp. free monoid) generated by the elements $\bomega_{a,b}$ with $0<b-a\le n$; it is convenient to set $\bomega_{a,a}=\bomega_{a,a+n+1}=\bold 1$ for $a\in\mathbb Z$.
 Given $r\ge 1$ and a multisegment $\bos= \{[c_1,d_1],\cdots,[c_r,d_r]\}\in{\rm mult}(\mathbb I_n)$, set $$\bomega_\bos=\bomega_{c_1,d_1}\cdots\bomega_{c_r,d_r}\in\cal I_n^+.$$ The simple objects of $\mathscr F_n$ are indexed by the elements of $\cal I_n^+$. Given $\bomega\in\cal I_n^+$, we let $V(\bomega)$ be the corresponding irreducible $\widehat\bu_n$-module. The objects of $\mathscr F_n$ are precisely the finite-dimensional representations whose Jordan--H\"older constituents are of the form $V(\bomega_\bos)$, $\bos\in{\rm mult}(\mathbb I_n)$.
\\\\
It was proved in \cite{HL10} that $\mathscr F_n$ is a monoidal tensor category. It follows from \cite{FR99} that the corresponding Grothendieck ring $\cal K(\mathscr F_n)$ is a polynomial ring (with coefficients in $\mathbb Z$) in the classes of the fundamental modules $[V(\bomega_{a,b})]$ with $0<b-a\le n$. \\\\ Given an object $V$ of $\mathscr F_n$ we let $\chi_q(V)\in\mathbb Z[\cal I_n]$ be the $q$-character (see \cite{FR99}) of $V$. Here $\mathbb Z[\cal I_n]$ is the group ring of $\cal I_n$, and we have
$$V=\bigoplus_{\bomegas\in\cal I_n} V_\bomegas, \qquad \chi_q(V)=\sum_{\bomegas\in\cal I_n}\dim V_\bomegas\ \bomega,$$
where $V_\bomegas$ is the $\ell$-weight space of $V$ of $\ell$-weight $\bomega$. It was proved in \cite{FR99} that the map $[V]\mapsto\chi_q(V)$ is an injective homomorphism of rings.

\subsection{Standard modules}\label{std} Given $\bos=\{[c_1,d_1],\cdots,[c_r,d_r]\}\in{\rm mult}(\mathbb I_n)$,
we define the standard module
$$W(\bomega_\bos)=V(\bomega_{c_1,d_1})\otimes\cdots\otimes V(\bomega_{c_r,d_r}).$$
It follows from \cite{CP94,ch01} (see also \cite[Corollary 3.2]{BC25b} for this formulation) that
$W(\bomega_{\bos})$ has simple head isomorphic to $V(\bomega_{\bos})$.
\\\\
Let
$c_{(1)}\ge\cdots\ge c_{(r)}$ be an enumeration of the multiset $\{c_s: 1\le s\le r\}$ in decreasing order and define
$\sigma\in\Sigma_r$ recursively by
$$\sigma(r)=\min\{k: c_{(k)}\le d_r\},\qquad
\sigma(p)=\min\{k\notin\{\sigma(p+1),\dots,\sigma(r)\}: c_{(k)}\le d_p\},\ \ 1\le p<r.$$
Set $\bos^\sharp=\{[c_{(\sigma(1))},d_1],\cdots,[c_{(\sigma(r))},d_r]\}$, where intervals of
the form $[c,c]$ are omitted.
No two segments of $\bos^\sharp$ are connected and hence it follows from \cite{ch01} that $$W(\bomega_{\bos^\sharp})\cong V(\bomega_{\bos^\sharp}).$$
The following was proved in \cite[Theorem 2(ii), Proposition 1.4.4]{BC25b}.  In  Section \ref{socstd} we give a different
proof, using normal sequences.
\begin{prop}\label{socW} Let $\bos\in{\rm mult}(\mathbb I_n)$ be such that $n\ge n(\bos)$. Then $W(\bomega_\bos)$ has
simple socle isomorphic to $V(\bomega_{\bos^\sharp})$.
\end{prop}

\subsection{The subcategories $\mathscr F_n(\bou_c)$}\label{bou} Given
$\bor=([c_1,d_1],\cdots,[c_r,d_r])\in{\rm mult}(\mathbb I_n)$ let $\mathscr F_n(\bor)$ be
the full subcategory of $\mathscr F_n$ consisting of objects whose Jordan--H\"older constituents are of
the form $V(\bomega)$, for $\bomega$ in the submonoid of $\cal I_n^+$ generated by $\bomega_{c_i,d_j}$,
$1\leq i,j\leq r$. It was proved in \cite[Proposition 5.1]{BC25b} that this is a monoidal subcategory of $\mathscr F_n$.
\\\\
We shall be interested in the following special case. Given $c\in\mathbb Z$, set
$\bou_c=([c,c+1],\cdots,[c+n-1,c+n])$. Then the simple objects of $\mathscr F_n(\bou_c)$ are the modules
$V(\bomega_\bos)$ with all segments of $\bos$ contained in $[c,c+n]$. In particular, if
$\boT_1,\cdots,\boT_k\in{\rm mult}(\mathbb I_n)$ and $n\ge n(\boT_1+\cdots+\boT_k)$, then taking
$c=\min\{a:[a,b]\in\boT_1+\cdots+\boT_k\}$, every composition factor of
$V(\bomega_{\boT_1})\otimes\cdots\otimes V(\bomega_{\boT_k})$ is an object of $\mathscr F_n(\bou_c)$.

\subsection{The involution $\Phi$}
For all $a,b\in\mathbb Z$ with $a<b$  set
\begin{gather}\label{omegaab}\bomega(a,b)=\bomega_{a,a+1}\bomega_{a+1,a+2}\cdots\bomega_{b-1,b}.\end{gather}
\begin{lem}\label{Phi}
    The assignment $$[V(\bomega_{a,b})]\mapsto [V(\bomega(a,b))],\ \ 0< b-a\le n,$$
defines an involution $\Phi: \cal K(\mathscr F_n)\to \cal K(\mathscr F_n)$.
\end{lem}
\begin{pf} Since $\cal K(\mathscr F_n)$ is a polynomial algebra in the elements $[V(\bomega_{a,b})]$, $0<b-a\le n$, it follows that $\Phi$ defines a homomorphism of rings.
    The following identity in $\cal K(\mathscr F_n)$ is a special case of \cite[Proposition 7.1]{BC25a}:
\begin{equation}\label{det1}[V(\bomega(a,b))]=\sum_{p=1}^{b-a} (-1)^{p-1}[V(\bomega_{a, a+p})][V(\bomega(a+p,b))].\end{equation}
Define $(b-a+1)\times (b-a+1)$ upper triangular matrices $E, H$ with rows and columns indexed by $\{a,a+1,\cdots, b\}$ whose $(i,j)$-th entries are the following elements of $\cal K(\mathscr F_n)$:
 $$E_{i,j}=(-1)^{j-i}[V(\bomega_{i,j})],\ \
H_{i, j}=[V(\bomega(i,j))],\ \ a\le i\le j\le b. $$ The preceding identity can then be reformulated as $EH=I$. It follows that $HE=I$, and the $(a,b)$
entry of $HE=I$ yields
\begin{equation}\label{det2}[V(\bomega_{a,b})]=\sum_{p=1}^{b-a}(-1)^{p-1}[V(\bomega(a, a+p))][V(\bomega_{a+p,b})].\end{equation}
To prove that $\Phi$ is an involution we proceed by induction on $b-a$.   Note that induction begins since  $V(\bomega_{a,a+1})=V(\bomega(a,a+1))$.
Applying $\Phi$ to both sides of \eqref{det1} and using induction on $b-a$ gives
\begin{gather*}\Phi([V(\bomega(a,b))])=\sum_{p=1}^{b-a}(-1)^{p-1}[V(\bomega(a,a+p))]\,\Phi([V(\bomega(a+p,b))])\\
=\sum_{p=1}^{b-a}(-1)^{p-1}[V(\bomega(a,a+p))]\,[V(\bomega_{a+p,b})]=[V(\bomega_{a,b})],\end{gather*} where the last equality is just \eqref{det2}, and the inductive step is established.
\end{pf}
\subsection{The co-standard module $M(\bos)$ }\label{zdual1} Given
$\bos=\{[c_1,d_1],\cdots,[c_r,d_r]\}\in{\rm mult}(\mathbb I_n)$, set
\begin{gather}
\bomega(c,c)=1,\ \ \bomega(\bos)= \bomega(c_1,d_1)\cdots\bomega(c_r,d_r),\\
M(\bos)= V(\bomega(c_1,d_1))\otimes\cdots\otimes V(\bomega(c_r,d_r)).\end{gather} Notice that if $\sigma\in\Sigma_r$ is such that $c_{\sigma(s)}\le d_s$ for all $1\le s\le r$, then
\begin{gather}\label{sig} \bomega(\bos)=\bomega(c_{\sigma(1)},d_1)\cdots\bomega(c_{\sigma(r)}, d_r).
\end{gather}
Since $\Phi$ is a ring homomorphism, we have $\Phi([W(\bomega_{\bos})])= [M(\bos)]$.
\subsection{The Zelevinsky dual}\label{zdual} Retain the notation of this section and, setting $b_1=d_1$, define
a sequence $a_1 > a_2 > \cdots > a_m$ recursively by
$$a_1= c_1, \ \ a_{s+1} = \max\{c_k : d_k = b_1 - s \ \ \text{and}\ \  c_k < a_s\},\ \ s\ge 1,$$
where $m \geq 1$ is maximal such that the set
$\{c_k : d_k = b_1-s \ \ \text{and}\ \  c_k < a_s\}$ is non-empty for each $1 \leq s \leq m-1$. In particular we have $a_{s+1}<b_1-s$ for all $1\le s\le m-1$.
Set
$$ \bos[1]= (\bos \setminus \{[a_1,b_1],\ldots,[a_m,b_1-m+1]\}) \cup \{[a_1,b_1-1],\ldots,[a_m,b_1-m]\},$$
where intervals of the form $[c,c]$ are omitted. Then $\bos[1]\in{\rm mult}(\mathbb I_n)$, and
the \emph{Zelevinsky dual} $\bos^\dagger$ is defined recursively by
\begin{equation}\label{lead}\bos^\dagger = ([b_1-m, b_1]) + (\bos[1])^\dagger,\end{equation}
with base case $\emptyset^\dagger = \emptyset$. 
When needed, we write $\bos^\dagger$ in the order of Section \ref{ms}. Note that
$\bomega_{\bos^\dagger}$, and hence $V(\bomega_{\bos^\dagger})$, depends only on $\bos^\dagger$ as a
multiset. 
\subsubsection{}
The following result is used repeatedly without mention.
\begin{lem}\label{daggerinv}
    For all $\bos\in\mult\In$ we have $$\bomega(\bos)=\bomega(\bos^\dagger),\ \ n(\bos)=n(\bos^\dagger),\ \ |\bos|=|\bos^\dagger|.$$
\end{lem}
\begin{pf} All statements are proved  by
induction on $|\bos|$, induction clearly beginning at $|\bos|=0$. Since
$|\bos[1]|=|\bos|-m<|\bos|$, we have
\begin{gather*}\bomega(\bos^\dagger)=\bomega(\bos[1]^\dagger)\,\bomega(b_1-m,b_1)=\bomega(\bos[1])\,\bomega(b_1-m,b_1)\\
=\bomega(\bos)\,\bomega(b_1-m,b_1)\prod_{p=1}^m\bomega(a_p,b_1-p)\,\bomega(a_p,b_1-p+1)^{-1}=\bomega(\bos).
\end{gather*}
Here the first equality is by definition of $\bos^\dagger$, the second is the inductive hypothesis, the third is
from the definition of $\bos[1]$, and the fourth follows by noting that $\bomega(c,d)=\bomega(c,d-1)\bomega_{d-1,d}$
for all $c<d$, so that the product over $p$ is $\prod_{p=1}^m\bomega_{b_1-p,b_1-p+1}^{-1}=\bomega(b_1-m,b_1)^{-1}$.  The proof of the other two statements is similar and we omit the details.
\end{pf}
\subsubsection{} For $\bos=\{[c_1,d_1],\cdots,[c_r,d_r]\}\in{\rm mult}(\mathbb I_n)$ let $\bos^\circ$ be the multisegment
of unit intervals given by
$$\mu_{\bos^\circ}([j,j+1])=\#\{s:\ c_s\le j<d_s\},\qquad j\in\mathbb Z.$$
Equivalently, $\bos^\circ$ is obtained by cutting every segment of $\bos$ into unit intervals and so
$$\bomega(\bos)=\bomega_{\bos^\circ}.$$ 
 For example, if $\bos=\{[0,2],[1,3]\}$ then $\bos^\circ=\{[0,1],[1,2],[1,2],[2,3]\}$ and $\bos^\sharp=\{[0,3],[1,2]\}$;
 one checks using the algorithm above that $\{[0,3],[1,2]\}^\dagger=\{[0,1],[1,2],[1,2],[2,3]\}$, i.e.
 $(\bos^\sharp)^\dagger=\bos^\circ$ (see Theorem \ref{zdualsoc}(ii)).
 
\subsection{} We state the main result of the paper.

\begin{thm}\label{zdualsoc}
Let $\bos\in{\rm mult}(\mathbb I_n)$. Then:
\begin{enumerit}
\item[(i)] $\operatorname{hd}M(\bos)\cong V(\bomega(\bos))\cong M(\bos^\sharp)$.
\item[(ii)] If $n\ge n(\bos)$, we have  $$\soc M(\bos)\cong V(\bomega_{\bos^\dagger}),\ \ \Phi([V(\bomega_\bos)])=[V(\bomega_{\bos^\dagger})], \ \ \Phi([V(\bomega_{\bos^\sharp})])=[V(\bomega(\bos))].$$

\end{enumerit}
\end{thm}
\noindent Using Proposition \ref{socW} and part (i) of the theorem we see that part (ii) shows, in particular, that the involution $\Phi$ exchanges the head
(resp. socle) of the standard module $W(\bomega_\bos)$ with the socle (resp. head) of the co-standard module
$M(\bos)$. In particular, $\bos^\sharp=(\bos^\circ)^\dagger$ depends only on the multiplicities
$\#\{s: c_s\le j<d_s\}$, $j\in\mathbb Z$. \\\\
The third equality in part (ii) is a simple consequence of the other statements of the theorem. Since $\bomega(\bos)=\bomega_{\bos^\circ}$ (see Section \ref{zdual}) it
suffices to prove that $(\bos^\sharp)^\dagger=\bos^\circ$. For this we note that $n(\bos^\sharp)= n(\bos)\le n$, since the
largest right end point and the smallest left end point of $\bos$ occur in segments of $\bos^\sharp$ which are not
omitted; hence the theorem applies to $V(\bomega_{\bos^\sharp})$. Since $W(\bomega_{\bos^\sharp})\cong V(\bomega_{\bos^\sharp})$
by Section \ref{std}, parts (ii) and (i) of Theorem \ref{zdualsoc} give
$$[V(\bomega_{(\bos^\sharp)^\dagger})]=\Phi([W(\bomega_{\bos^\sharp})])=[M(\bos^\sharp)]=[V(\bomega(\bos))]
=[V(\bomega_{\bos^\circ})].$$
Hence $\bomega_{(\bos^\sharp)^\dagger}=\bomega_{\bos^\circ}$, and since $\cal I_n^+$ is a free monoid this
gives $(\bos^\sharp)^\dagger=\bos^\circ$.\\\\ 
The proof of the rest of the  theorem is given in Sections \ref{head}--\ref{phi}. 
\subsection{Multiplicities in co-standard modules}\label{multdual} We note the following immediate consequence of Theorem \ref{zdualsoc}. For $\bos=\{[c_1,d_1],\cdots,[c_r,d_r]\}\in{\rm mult}(\mathbb I_n)$ let
$\lambda_\bos=\sum_{j=1}^r\omega_{d_j-c_j}$, where $\omega_1,\cdots,\omega_n$ are the fundamental weights of
$\lie{sl}_{n+1}$. It is well known (see \cite{CP01}) that the $\lambda_\bos$-weight space of $W(\bomega_\bos)$ is
one-dimensional, and hence
\begin{equation}\label{Wexp}[W(\bomega_\bos)]=[V(\bomega_\bos)]+\sum_{\boT}m_{\bos,\boT}\,[V(\bomega_{\boT})],\qquad m_{\bos,\boT}\in\mathbb Z_+,\end{equation}
where the sum is over $\boT\in{\rm mult}(\mathbb I_n)$ with $\lambda_{\boT}<\lambda_\bos$ with respect to the usual
order on dominant weights of $\lie{sl}_{n+1}$. Here $m_{\bos,\boT}$ is the multiplicity of $V(\bomega_\boT)$
as a composition factor of $W(\bomega_\bos)$, and $m_{\bos,\bos}=1$. Applying $\Phi$ to both sides of the equality in \eqref{Wexp} and using $\Phi([W(\bomega_\bos)])=[M(\bos)]$ and Theorem \ref{zdualsoc}(ii) we get  $$[M(\bos)]=[V(\bomega_{\bos^\dagger})]+\sum_{\boT}m_{\bos,\boT}\,[V(\bomega_{\boT^\dagger})].$$
In particular, $[M(\bos):V(\bomega_{\boT^\dagger})]=[W(\bomega_\bos):V(\bomega_\boT)]$ for all $\boT$,
and every composition factor of $M(\bos)$ is of the form $V(\bomega_{\boT^\dagger})$ with
$\lambda_\boT\le\lambda_\bos$.
\subsection{Real and prime modules} Recall that a simple module $V$ is \emph{prime} if it is not isomorphic to a tensor product of two
non-trivial modules, and \emph{real} if $V\otimes V$ is irreducible; otherwise $V$ is \emph{imaginary}.
\begin{prop}\label{realprime} Let $\bos\in{\rm mult}(\mathbb I_n)$ with $n\ge n(\bos)$. Then
$V(\bomega_\bos)$ is real (resp. prime) if and only if $V(\bomega_{\bos^\dagger})$ is real (resp. prime).
\end{prop}
\begin{pf} Since $(\bos^\dagger)^\dagger=\bos$ and $n(\bos^\dagger)=n(\bos)$ (Lemma \ref{daggerinv}), it suffices to prove the forward direction. Suppose that $V(\bomega_\bos)\cong V(\bomega_{\bos_1})\otimes V(\bomega_{\bos_2})$ for some
$\bos_1,\bos_2\in{\rm mult}(\mathbb I_n)$ with $\bos=\bos_1+\bos_2$. Then $n\ge n(\bos_i)$ and so applying $\Phi$ and using
Theorem \ref{zdualsoc}(ii) gives
$$[V(\bomega_{\bos^\dagger})]=[V(\bomega_{\bos_1^\dagger})][V(\bomega_{\bos_2^\dagger})],\ \ {\rm i.e.}\ \ V(\bomega_{\bos^\dagger})\cong V(\bomega_{\bos_1^\dagger})\otimes V(\bomega_{\bos_2^\dagger}).$$
 Taking $\bos_1,\bos_2$
non-empty (so that $\bos_1^\dagger,\bos_2^\dagger$ are non-empty) shows that if $V(\bomega_\bos)$ is not prime then neither is
$V(\bomega_{\bos^\dagger})$; taking $\bos_1=\bos_2=\bos$ (note $n(\bos+\bos)=n(\bos)$) shows that if $V(\bomega_\bos)$ is real then
so is $V(\bomega_{\bos^\dagger})$.
\end{pf}
In \cite{BC23} and \cite{BM25} the authors construct large families of imaginary modules as
subquotients of tensor products of snake modules. The examples in these papers always correspond to multisegments with an even number of elements. Proposition \ref{realprime} applied to these gives further
imaginary modules which may not arise from these constructions. We give one example.\\\\
{\bf Example.} \label{primeex} Let
\begin{gather*} \bos=\{[-4,0],[-2,3],[0,2],[0,5],[2,4],[4,7]\},\\
\bos^\dagger=\{[-4,-3],[-3,-1],[-2,1],[-1,2],[0,1],[1,3],[1,3],[2,5],[3,4],[4,6],[6,7]\}.\end{gather*}
By \cite[Theorem 3]{BC23}, for $n=6$ the module $V(\bomega_\bos)$ occurs in the Jordan--H\"older series of
$V(\bomega_{0,3}\bomega_{2,5}\bomega_{4,7})\otimes V(\bomega_{-4,0}\bomega_{-2,2}\bomega_{0,4})$ and is imaginary; restricting to the
subalgebra associated to the first six nodes shows that it is imaginary for all $n\ge6$. Since $n(\bos)=11$,
Proposition \ref{realprime} shows that $V(\bomega_{\bos^\dagger})$ is imaginary for all $n\ge11$; it remains open whether it is
imaginary for $6\le n\le 10$.
\\\\
Since $\bos^\dagger$ has an odd number of elements it does not arise directly from the work of \cite{BC23} or \cite{BM25}. Further neither $\bos$ nor $\bos^\dagger$ is regular: $\bos$ has a repeated left end point, and $\bos^\dagger$ has the segment $[1,3]$
twice. In particular neither of these modules is covered by the criterion of \cite{LM18}, which applies to regular
multisegments.
\\\\
Finally, $V(\bomega_\bos)$ is prime for $n\ge6$, and hence by Proposition \ref{realprime} so is $V(\bomega_{\bos^\dagger})$ for
$n\ge11$. We only sketch the argument. A module whose multisegment has at most three segments is real (this
uses \cite{BC25a}). Hence if $V(\bomega_\bos)$ were not prime, it would be a tensor product
$V(\bomega_{\bos_1})\otimes V(\bomega_{\bos_2})$ with $\bos_2$ consisting of one or two segments. In each of these $21$ cases
one finds, by a computation with $q$-characters, an $\ell$-weight of $V(\bomega_{\bos_1})\otimes V(\bomega_{\bos_2})$ which is not
an $\ell$-weight of $V(\bomega_\bos)$, so the tensor product is not simple.

\subsection{Alternating Snake Modules and  a determinantal formula}\label{snakesdef} We end this section with a discussion of a particular family of modules: the alternating snake modules, which were introduced in \cite{BC25a}. These modules have a particularly nice character formula which is an analog of the classical Jacobi--Trudi formula and motivates the discussion in the next section.

\subsubsection{} Say that $\bos\in{\rm mult}(\mathbb I_n)$ is an alternating snake if
$\mu_\bos([c,d])\le 1$ for all $[c,d]\in\mathbb I_n$, and if $\bos$ can be written as an ordered tuple $([i_1,j_1],\cdots, [i_r, j_r])$ for some $r\ge 1$, and if the following hold: 
\begin{enumerit}
\item[(i)] for $1\le s\le r-1$ there exists $\epsilon_s\in\{0,1\}$ such that
$i_{s+\epsilon_s}<i_{s+1-\epsilon_s}$ and $j_{s+\epsilon_s}<j_{s+1-\epsilon_s}$,
\item[(ii)] if $1\le\ell\le r-2$ and $\epsilon_\ell\ne\epsilon_{\ell+1}$, then $[i_k,j_k]$ and $[i_m,j_m]$ do not overlap for all
$1\le k\le\ell$ and $\ell+2\le m\le r$.
\end{enumerit}
\begin{rem}\label{snakes} If $\epsilon_s=\epsilon_{1}$ for all $1\le s\le r-1$  then we say that $\bos$ is an increasing snake if $\epsilon_1=0$ and decreasing if $\epsilon_1=1$. The definition and study of snakes was begun in \cite{MY12}.

\end{rem}

\subsubsection{The determinantal formula}
Given $\bos=([i_1,j_1], \cdots, [i_r,j_r])$  an alternating snake, we define an $r\times r$ matrix $A(\bos)$ with entries in $\cal K(\mathscr F_n)$ recursively as follows. \\\\
If $r=1$ set $A(\bos)=([V(\bomega_{i_1,j_1})])$. Assume that $r\ge2$ and that $A(\bos')$ has been defined for alternating
snakes with $r-1$ segments.  In particular, it is defined for $\bos(1,r)=([i_2,j_2],\cdots, [i_r,j_r])$ which is also an alternating snake. Let $r_1\le r$ be maximal such that $\epsilon_1=\cdots=\epsilon_{r_1-1}$. Then
\begin{itemize}
\item $A(\bos)_{p,\ell}=A(\bos(1,r))_{p-1,\ell-1}$ for $2\le p,\ell\le r$;
\item $A(\bos)_{1,\ell}=[V(\bomega_{i_1,j_\ell})]$ if $[i_1,j_\ell]\in\mathbb I_n$, $([i_s,j_s],[i_{s+1},j_{s+1}])$ is connected for all
$1\le s\le\ell-1$ and either
$$\ell\le r_1,\qquad\text{or}\qquad r_1<\ell,\ \ \epsilon_1=1\ \ \text{and}\ \ \epsilon_{r_1}=\cdots=\epsilon_{\ell-1}=0,$$ and in all other cases $A(\bos)_{1,\ell}=0$.
\item $A(\bos)_{\ell,1}=[V(\bomega_{i_\ell,j_1})]$ if  $[i_\ell,j_1]\in\mathbb I_n$, $([i_s,j_s],[i_{s+1},j_{s+1}])$ is connected for all
$1\le s\le\ell-1$ and either
$$\ell\le r_1,\qquad\text{or}\qquad r_1<\ell,\ \ \epsilon_1=0\ \ \text{and}\ \ \epsilon_{r_1}=\cdots=\epsilon_{\ell-1}=1$$ and in all other cases $A(\bos)_{\ell,1}=0$.
\end{itemize}
In particular, if $\bos$ is a snake then $r_1=r$ and $\det A(\bos)=\det(e_{i_p,j_\ell})_{1\le p,\ell\le r}$, where we understand that
$e_{i_p,j_\ell}= [V(\bomega_{i_p, j_\ell})]$ if $[i_p,j_\ell]\in\mathbb I_n$ and zero otherwise.\\\\ Finally, following \cite[Section 2.8.3]{BC25a}, we say that $\bos$ is stable if for all $2\le p\le r-1$ the segment
$[i_{p+1},j_{p+1}]$ does not contain the segment $[i_{p-1},j_{p-1}]$.
\begin{thm*}[{\cite[Theorem 3]{BC25a}}]\label{altdet} If $\bos$ is a stable alternating snake, then
$[V(\bomega_\bos)]=\det A(\bos)$.
\end{thm*}

\section{An interval analogue of the ring of symmetric functions}\label{lambdasec}

The results of Section \ref{sec:main} suggest that $\cal K(\mathscr F_n)$ should be viewed as an analogue of the
ring of symmetric functions: the involution $\Phi$ of Lemma \ref{Phi} interchanges the classes of
the fundamental modules $V(\bomega_{a,b})$ and of the modules $V(\bomega(a,b))$, in the same way that the
classical involution $\omega$ interchanges the elementary and the complete symmetric functions, and
the determinantal formula of \cite[Theorem 3]{BC25a} for alternating snake modules (see Section
\ref{snakesdef}) is an analogue of the Jacobi--Trudi formula. In this section we make this precise
by constructing a ring ${\rm Sym}(\mathbb I_n)$ of symmetric functions in which the generators depend
on intervals rather than on positive integers, together with an isomorphism
${\rm Sym}(\mathbb I_n)\cong\cal K(\mathscr F_n)$.
\subsection{} Recall that $\mathbb I_n^+$ is the subset of $\mathbb I_n$ consisting of intervals $[a,b]$ with $a<b$. Denote by $\mathbb Z[\mathbb I^+_n]$ the polynomial ring over $\mathbb Z$ in variables $x_{a,b}$ with $[a,b]\in\mathbb I_n^+$.\\\\
For $[a,b]\in\mathbb I_n^+$ let ${\rm Fun}_{a,b}$ be the set of functions $$f:\{a+1,\cdots,b\}\to\mathbb Z,\ \ d-1<f(d)\le d+n,\ \ a<d\le b.$$ Let $\mathbb P_{a,b}$ (resp. $\mathbb Q_{a,b}$) be the subset of ${\rm Fun}_{a,b}$ consisting of weakly decreasing (resp. strictly increasing) functions.
Set $e_{a,a}=h_{a,a}=1$ and,
for $[a,b]\in\mathbb I_n^+$, set \begin{gather}\label{eab}
e_{a,b}=\sum_{p\in\mathbb P_{a,b}}x_{a,p(a+1)}x_{a+1,p(a+2)}\cdots x_{b-1,p(b)},\\
h_{a,b}=\sum_{q\in \mathbb Q_{a,b}}x_{a,q(a+1)}x_{a+1,q(a+2)}\cdots x_{b-1,q(b)}.
\end{gather}
In particular $e_{a,a+n+1}=x_{a,a+n+1}x_{a+1,a+n+1}\cdots x_{a+n,a+n+1}$.
\begin{prop}\label{ehrel}
    The elements $e_{a,b}$, $0<b-a\le n$, are algebraically independent, and for $[a,b]\in\mathbb I_n^+$ with $b-a\le n$ we have $$\sum_{k=0}^{b-a}(-1)^k\,e_{a,a+k}\,h_{a+k,b}=0.$$
\end{prop}
The proof of the proposition is given in two steps.
\subsubsection{} Let $\mathbb Z[\cal I_n]$ be the group ring of the free abelian group
$\cal I_n$. Recall from Section \ref{ms} that we have an injective homomorphism of rings $\chi_q:\cal K(\mathscr F_n)\to\mathbb Z[\cal I_n]$. (In the notation of \cite{FR99} the element
$\bomega_{a,b}$ corresponds to $Y_{b-a,q^{a+b}}$.)\\\\ Let
$\phi:\mathbb Z[\mathbb I_n^+]\to\mathbb Z[\cal I_n]$ be the homomorphism of rings defined by extending the assignment
$x_{a,b}\mapsto\bomega_{a,b}\bomega_{a+1,b}^{-1}$, $ [a,b]\in\mathbb I_n^+$.
Clearly, $$\phi(x_{a,b}x_{a+1,b}\cdots x_{b-1,b})=\bomega_{a,b}.$$
The following is a reformulation of the results of \cite{MY12}.
\begin{lem}\label{qchar} For $[a,b]\in\mathbb I_n^+$ we have
$$\phi(e_{a,b})=\chi_q(V(\bomega_{a,b})),\qquad \phi(h_{a,b})=\chi_q(V(\bomega(a,b))).$$
\end{lem}
\begin{pf} The first statement of the lemma is immediate from the description of $q$-characters given in \cite{MY12}. For the second, we note that in their language the module
$V(\bomega(a,b))$ is associated to the snake
$([a,a+1],[a+1,a+2],\cdots,[b-1,b])$ of unit segments. According to their result the $q$-character of this module is given by the tuples $(p_0,\cdots, p_{b-a-1})$ with $p_s\in\mathbb P_{a+s,a+s+1}$, $0\le s\le b-a-1$, such that $p_s(a+s+1)<p_{s+1}(a+s+2)$. It is now trivial to see that such tuples are in bijective correspondence with the elements of $\mathbb Q_{a,b}$.
\end{pf}
\medskip
The elements $\chi_q(V(\bomega_{a,b}))$
are algebraically independent since $\chi_q$ is injective and $\cal K(\mathscr F_n)$ is the polynomial ring over
$\mathbb Z$ in the classes $[V(\bomega_{a,b})]$ with $0<b-a\le n$. It follows immediately from the lemma that the elements $e_{a,b}$, $0<b-a\le n$,
 are also algebraically independent, and the first assertion of the proposition is proved.

\subsubsection{} Suppose that $f_1, f_2\in{\rm Fun}_{a,b}$. We claim that $$\phi(x_{a,f_1(a+1)}x_{a+1,f_1(a+2)}\cdots x_{b-1,f_1(b)})=\phi(x_{a,f_2(a+1)}x_{a+1,f_2(a+2)}\cdots x_{b-1,f_2(b)})\iff f_1=f_2.$$
The reverse direction is obvious. For the forward direction we are given that
$$\prod_{c=a}^{b-1}\bomega_{c, f_1(c+1)}\bomega_{c+1,f_1(c+1)}^{-1}=\prod_{c=a}^{b-1}\bomega_{c, f_2(c+1)}\bomega_{c+1,f_2(c+1)}^{-1}.$$ Since a term of the form $\bomega_{a,j}$ occurs on the left hand side if and only if $j=f_1(a+1)<a+n+1$, and on the right hand side if and only if $j=f_2(a+1)<a+n+1$, it follows that $f_1(a+1)=f_2(a+1)$. Hence we can cancel the appropriate terms on both sides to get an equality where the product begins with $a+1$. An obvious iteration gives the claim.\\\\
Suppose now that $f_1,\cdots, f_k$ are distinct elements of ${\rm Fun}_{a,b}$ and $c_1,\cdots, c_k\in\mathbb Z$ are such that
$$\phi\left( \sum_{s=1}^kc_sx_{a,f_s(a+1)}x_{a+1,f_s(a+2)}\cdots x_{b-1,f_s(b)}\right)=0.$$
By the preceding discussion the elements $\phi(x_{a,f_s(a+1)}x_{a+1,f_s(a+2)}\cdots x_{b-1,f_s(b)})$ are distinct elements of $\cal I_n$, and hence linearly independent in $\mathbb Z[\cal I_n]$; therefore $c_s=0$ for all $s$. In other words, $\phi$ is injective on the span of the monomials in the set \begin{equation}\label{basis} S_{a,b}= \{x_{a,f(a+1)}x_{a+1,f(a+2)}\cdots x_{b-1,f(b)}: f\in{\rm Fun}_{a,b}\}, \end{equation} and the elements of $S_{a,b}$ are linearly independent.
A pair $(p',q')\in\mathbb P_{a,a+k}\times\mathbb Q_{a+k,b}$ defines an element of
${\rm Fun}_{a,b}$ in an obvious way. Hence the elements $$e_{a,a+k}h_{a+k,b},\ \ 0\le k\le b-a,\ \ {\rm and\ so}\ \  \sum_{k=0}^{b-a}(-1)^k\,e_{a,a+k}\,h_{a+k,b},$$ are in the span of $S_{a,b}$. Lemma \ref{qchar} and \eqref{det1} show that
$$\phi\Big(\sum_{k=0}^{b-a}(-1)^k\,e_{a,a+k}\,h_{a+k,b}\Big)=0.$$
Hence $\sum_{k=0}^{b-a}(-1)^k\,e_{a,a+k}\,h_{a+k,b}=0$, as needed. The proof of Proposition \ref{ehrel} is complete.

\subsection{The ring $\Sym(\In)$}\label{Phisym} Let ${\rm Sym}(\mathbb I_n)$ be the subring of $\mathbb Z[\mathbb I_n^+]$
generated by the elements $e_{a,b}$, $1\le b-a\le n$. Proposition \ref{ehrel} shows that it is a
polynomial ring over $\mathbb Z$ in these elements, and that it contains the
elements $h_{a,b}$, $1\le b-a\le n$. \\\\
Lemma \ref{qchar} shows that the map
$$\Phi_n=\chi_q^{-1}\circ\phi:{\rm Sym}(\mathbb I_n)\to\cal K(\mathscr F_n)$$
is an isomorphism of rings with
$$\Phi_n(e_{a,b})=[V(\bomega_{a,b})],\qquad \Phi_n(h_{a,b})=[V(\bomega(a,b))],\qquad 1\le b-a\le n.$$
Recalling the involution $\Phi$ of Lemma \ref{Phi}, we note that $\Phi_n^{-1}\circ\Phi\circ\Phi_n$ defines an
involution of ${\rm Sym}(\mathbb I_n)$ with $e_{a,b}\mapsto h_{a,b}$, which by abuse of
notation we continue to denote by $\Phi$. In particular, the elements $h_{a,b}$,
$1\le b-a\le n$, are also algebraically independent and generate ${\rm Sym}(\mathbb I_n)$. The following is now immediate.
\begin{prop} The elements
$${\rm sch}(\bos)=\Phi_n^{-1}([V(\bomega_{\bos})]),\qquad \bos\in\mult(\In),$$
form a $\mathbb Z$-basis of ${\rm Sym}(\mathbb I_n)$. Moreover ${\rm sch}([a,b])=e_{a,b}$ and
${\rm sch}([a,a+1],\cdots,[b-1,b])=h_{a,b}$.\hfill\qedsymbol
\end{prop}

\subsubsection{Positivity} Let $\mathbb Z[\mathbb I_n^+]^\flat$ be the quotient of
$\mathbb Z[\mathbb I_n^+]$ by the ideal generated by the elements $e_{a,a+n+1}-1$, $a\in\mathbb Z$, where
$$e_{a,a+n+1}=x_{a,a+n+1}x_{a+1,a+n+1}\cdots x_{a+n,a+n+1}.$$ It
follows that $\mathbb Z[\mathbb I_n^+]^\flat$ is isomorphic to the Laurent polynomial ring in the
variables $x_{a,b}^{\pm1}$ with $1\le b-a\le n$. Since $\phi(e_{a,a+n+1})=\bomega_{a,a+n+1}=\bold 1$, we have an induced
homomorphism of rings
$$\bar\phi:\mathbb Z[\mathbb I_n^+]^\flat\to\mathbb Z[\cal I_n],\qquad
x_{a,b}\mapsto\bomega_{a,b}\bomega_{a+1,b}^{-1},\ \ 1\le b-a\le n,$$
where we recall that $\bomega_{a+1,a+1}=\bold 1$. In fact $\bar\phi$ is an isomorphism, with inverse given by
$\psi(\bomega_{a,b})=x_{a,b}x_{a+1,b}\cdots x_{b-1,b}$, $1\le b-a\le n$.
Finally, since $\phi$ is injective on ${\rm Sym}(\mathbb I_n)$ and factors through
$\mathbb Z[\mathbb I_n^+]^\flat$, it follows that we have an induced injective map ${\rm Sym}(\mathbb I_n)\to\mathbb Z[\mathbb I_n^+]^\flat$.

\begin{prop}\label{positive} Let $\bos$ be a multisegment. The image of ${\rm sch}(\bos)$ in
$\mathbb Z[\mathbb I_n^+]^\flat$ is
$$\sum_M\dim V(\bomega_\bos)_{\bar\phi(M)}\ M,$$
where $M$ varies over the Laurent monomials in the variables $x_{a,b}$, $1\le b-a\le n$. In
particular, it is a linear combination of Laurent monomials with non-negative integer
coefficients.
\end{prop}
\begin{pf} Write the image of ${\rm sch}(\bos)$ as $\sum_Ma_MM$. Applying $\bar\phi$ gives
$$\sum_Ma_M\,\bar\phi(M)=\phi({\rm sch}(\bos))=\chi_q(V(\bomega_\bos))
=\sum_{\bomega}\dim V(\bomega_\bos)_{\bomega}\ \bomega,$$
and since the elements $\bar\phi(M)$ are pairwise distinct, $a_M=\dim V(\bomega_\bos)_{\bar\phi(M)}$.
\end{pf}

\subsection{Connection to the classical literature}\label{classical} Let $\Lambda$ be the ring of symmetric functions
in infinitely many variables and $\omega$ its classical involution, $\omega(e_k)=h_k$. For $N\ge1$ let
$\Lambda_N=\mathbb Z[z_1,\dots,z_N]^{\Sigma_N}$ be the ring of symmetric polynomials in $N$ variables, and let $\rho_N:\Lambda\to\Lambda_N$ be
the restriction to $N$ variables; it is surjective, with kernel the ideal generated by the $e_k$, $k>N$, so that $\Lambda_N$ is a
quotient of $\Lambda$.
We study 
the relation between ${\rm Sym}(\mathbb I_n)$ and these rings.
\subsubsection{The homomorphism $\pi$} Consider the homomorphism
$$\pi:\mathbb Z[\mathbb I_n^+]\to\mathbb Z[z_1,\dots,z_{n+1}],\qquad x_{a,b}\mapsto z_{b-a}.$$ The following is trivially checked.
\begin{prop} For $0\le b-a\le n+1$ we have
$$\pi(e_{a,b})=e_{b-a}(z_1,\dots,z_{n+1}),\qquad \pi(h_{a,b})=h_{b-a}(z_1,\dots,z_{n+1}).$$
In particular $\pi(e_{a,a+n+1})=z_1\cdots z_{n+1}$ and $\pi({\rm Sym}(\mathbb I_n))\subset\Lambda_{n+1}$. \hfill\qedsymbol
\end{prop}
\subsubsection{The map $\kappa$} Recall from Section \ref{Phisym}  the involution $\Phi$ of ${\rm Sym}(\mathbb I_n)$ with $\Phi(e_{a,b})=h_{a,b}$. To compare it
with the involution $\omega$ of $\Lambda$ we lift $\pi$ to $\Lambda$. Since ${\rm Sym}(\mathbb I_n)$ is a polynomial ring in the elements
$e_{a,b}$, $1\le b-a\le n$, there is a unique homomorphism of rings
$$\kappa:{\rm Sym}(\mathbb I_n)\to\Lambda,\qquad \kappa(e_{a,b})=e_{b-a},\ \ 1\le b-a\le n.$$

\begin{lem}\label{kappa} We have $\kappa(h_{a,b})=h_{b-a}$ for $1\le b-a\le n$. Moreover
$$\rho_{n+1}\circ\kappa=\pi|_{\Sym(\In)}\qquad \kappa\circ\Phi=\omega\circ\kappa,$$
where $\Phi$ is the involution of ${\rm Sym}(\mathbb I_n)$. In other words, $\kappa$ intertwines the involutions $\Phi$ of
${\rm Sym}(\mathbb I_n)$ and $\omega$ of $\Lambda$.
\end{lem}
\begin{pf} Using  Proposition \ref{ehrel} and the result in the classical case, we have
$$\sum_{p=0}^{b-a}(-1)^pe_{a,a+p}h_{a+p,b}=0= \sum_{p=0}^{b-a}(-1)^p e_ph_{b-a-p},\ \ b-a>0. $$ An obvious induction on $b-a$ gives $\kappa(h_{a,b})=h_{b-a}$. The equalities hold since the homomorphisms  agree on the generators $e_{a,b}$.
\end{pf}
\noindent In ${\rm Sym}(\mathbb I_n)$, Theorem \ref{zdualsoc}(ii)
says that $\Phi({\rm sch}(\bos))={\rm sch}(\bos^\dagger)$ if $n\ge n(\bos)$; applying $\kappa$ and using Lemma \ref{kappa} we get, for
$n\ge n(\bos)$,
\begin{equation}\label{omegas}\omega(\kappa({\rm sch}(\bos)))=\kappa(\Phi({\rm sch}(\bos)))=\kappa({\rm sch}(\bos^\dagger)),\end{equation}
the analogue of the identity $\omega(s_\lambda)=s_{\lambda'}$. \noindent In particular $c_\lambda(\bos)\ge0$ whenever $\ell(\lambda)\le n+1$. If $\kappa({\rm sch}(\bos))$ does not depend on the choice of $n\ge n(\bos)$, then
taking $n\ge|\bos|-1$ shows that all the coefficients $c_\lambda(\bos)$ are multiplicities, so that $\kappa({\rm sch}(\bos))$ is Schur positive.

\begin{rem} It is necessary to pass from $\Lambda_{n+1}$ to $\Lambda$ here, since $\omega$ does not preserve the kernel of $\rho_{n+1}$.
The assignment $e_k\mapsto h_k$, $1\le k\le N$, does define an involution of $\Lambda_N$, but it does not map $s_\lambda$ to $s_{\lambda'}$,
even for partitions $\lambda$ in an $N\times N$ box, since for $k>N$ the element $h_k$ is not one of the generators of
$\Lambda_N$. For instance, if $N=2$ then $s_{(2,2)}=h_2^2-h_1h_3$ and $h_3=2h_1h_2-h_1^3$ in $\Lambda_2$, so this involution maps $h_3$
to $-h_3$ and $s_{(2,2)}$ to $e_2^2+e_1h_3\ne e_2^2=s_{(2,2)}$.
\end{rem}

\subsection{The elements $\kappa({\rm sch}(\bos))$} We now relate $\kappa({\rm sch}(\bos))$ to characters. Let $\varepsilon_1,\dots,\varepsilon_{n+1}$ be the weights of the natural
representation of $\lie{gl}_{n+1}$, so that $z_i$ has weight $\varepsilon_i$, and let $P$ be the $\mathbb Z$-span of these weights
modulo the relation $\sum_{i=1}^{n+1}\varepsilon_i=0$, i.e.\ the weight lattice of $\lie{sl}_{n+1}$. The assignment
$\bomega_{a,b}\mapsto \varepsilon_1+\cdots+\varepsilon_{b-a}$, $[a,b]\in\mathbb I_n$, defines a
homomorphism of group rings $$\mathbb Z[\cal I_n]\to\mathbb Z[P], \ \ \ \ \chi_q(V)\mapsto 
\ch_{\bu_n}V.$$  It follows from Proposition \ref{positive} and Lemma
\ref{kappa} that
$$\rho_{n+1}(\kappa({\rm sch}(\bos)))=\pi({\rm sch}(\bos))\equiv\ch_{\bu_n}V(\bomega_\bos)\pmod{z_1\cdots z_{n+1}-1}.$$
\noindent Note that $e_{a,b}$ and $h_{a,b}$ are homogeneous polynomials of degree $b-a$ in the variables $x_{c,d}$, so that the
grading of $\mathbb Z[\mathbb I_n^+]$ by degree induces a grading of ${\rm Sym}(\mathbb I_n)$ in which $e_{a,b}$ has degree $b-a$; since
$\kappa(e_{a,b})=e_{b-a}$, the map $\kappa$ preserves degrees.
\begin{lem}\label{homog} Let $\bos\in{\rm mult}(\mathbb I_n)$ with $n\ge n(\bos)$. Then ${\rm sch}(\bos)$ is homogeneous of degree
$|\bos|=\sum_s(d_s-c_s)$, and hence so is $\kappa({\rm sch}(\bos))$.
\end{lem}
\begin{pf} If $\boT=\{[c_{\sigma(1)},d_1],\cdots,[c_{\sigma(r)},d_r]\}$, where $\sigma\in\Sigma_r$ satisfies $c_{\sigma(p)}\le d_p$ for all $p$ (intervals
of the form $[c,c]$ being omitted), then it is clear that $$|\boT|=|\bos|,\ \ n(\boT)\le n(\bos).$$
By \cite[Theorem 1(ii)]{BC25b} all the dominant $\ell$-weights of $W(\bomega_\boT)$ are of this form, and hence so is the
highest $\ell$-weight of every composition factor of $W(\bomega_\boT)$; here we use that $n\ge n(\bos)$, which guarantees that no
segment of length $n+1$ (and hence no factor $\bomega_{a,a+n+1}=\bold 1$) occurs. Hence the expansions \eqref{Wexp} of the
classes $[W(\bomega_\boT)]$, with $\boT$ of this form, only involve classes $[V(\bomega_{\boT'})]$ with $\boT'$ of this form, and in
particular with $|\boT'|=|\bos|$. Since these expansions are unitriangular and finite, inverting them gives
$$[V(\bomega_\bos)]=\sum_{\boT\in\mult(\In)}b_{\bos,\boT}\,[W(\bomega_\boT)],\qquad b_{\bos,\boT}\ne 0\implies |\boT|=|\bos|.$$
The lemma follows by applying $\Phi_n^{-1}$ and noting that $\Phi_n^{-1}([W(\bomega_\boT)])$ is the product of the elements $e_{a,b}$,
$[a,b]$ a segment of $\boT$, and hence is homogeneous of degree $|\boT|=|\bos|$.
\end{pf}
\noindent For a partition $\lambda$ with $\ell(\lambda)\le n+1$ let $V(\lambda)$ be the simple $U_q(\lie{sl}_{n+1})$-module with highest weight
$\sum_i(\lambda_i-\lambda_{i+1})\omega_i$. If $n\ge n(\bos)$, then by Lemma \ref{homog} and the preceding discussion,
$$\rho_{n+1}\big(\kappa({\rm sch}(\bos))\big)=\pi({\rm sch}(\bos))=\sum_{\substack{\lambda\vdash|\bos|\\ \ell(\lambda)\le n+1}}[V(\bomega_\bos):V(\lambda)]\,s_\lambda(z_1,\dots,z_{n+1}).$$
Since the polynomials $s_\lambda(z_1,\dots,z_{n+1})$, $\ell(\lambda)\le n+1$, are linearly independent and $\rho_{n+1}(s_\lambda)=0$ if
$\ell(\lambda)>n+1$, it follows that $$\kappa({\rm sch}(\bos))=\sum_{\lambda\vdash|\bos|}c_\lambda(\bos)s_\lambda,\ \ 
c_\lambda(\bos)=[V(\bomega_\bos):V(\lambda)],\ \ \ell(\lambda)\le n+1.$$

\begin{rem}\label{conjpart} If $n\ge n(\bos)$, then by \eqref{omegas} we have $\kappa({\rm sch}(\bos^\dagger))=\omega(\kappa({\rm sch}(\bos)))$, and so $c_{\lambda'}(\bos^\dagger)=c_\lambda(\bos)$ for
all $\lambda$. Hence the decompositions of $V(\bomega_\bos)$ and $V(\bomega_{\bos^\dagger})$ into simple $U_q(\lie{sl}_{n+1})$-modules are
related by conjugation of partitions, $[V(\bomega_{\bos^\dagger}):V(\lambda')]=[V(\bomega_\bos):V(\lambda)]$, for all $\lambda$ with $\ell(\lambda)\le n+1$ and $\lambda_1\le n+1$. These
conditions cannot be dropped: for $k\ge n+2$ copies of $[0,1]$ the multisegment is self-dual, but $(k)$ occurs while $(1^k)$ has too many rows.
\end{rem}
\subsection{Snakes and skew diagrams}\label{skew} For snakes the analogy between $\bos\mapsto\bos^\dagger$ and the conjugation of
partitions is literal. Let $\lambda=(\lambda_1\ge\cdots\ge\lambda_\ell)$ and $\mu=(\mu_1\ge\cdots\ge\mu_\ell)$ be partitions with $\mu_i\le\lambda_i$ for all $i$,
and let $\lambda/\mu$ be the corresponding skew Young diagram; its transpose is $\lambda'/\mu'$. The associated snakes are the multisegments
$$\bor(\lambda/\mu)=\{[\mu_i-i+1,\lambda_i-i+1]:\ \mu_i<\lambda_i\},\qquad \boc(\lambda/\mu)=\{[j-\lambda'_j,\,j-\mu'_j]:\ \mu'_j<\lambda'_j\}.$$ As an example, the rectangle $\lambda=(i^r)$, $\mu=0$, gives $\bor(\lambda)=\bok_a(i,r)$ and $\boc(\lambda)=\bok_a(r,i)$ with $a=1-r$, where
$$\bok_a(i,r)=\{[a,a+i],[a+1,a+i+1],\cdots,[a+r-1,a+r-1+i]\}$$
is the multisegment of a Kirillov--Reshetikhin module.\\\\
Conversely, if $\bos=([a_1,b_1],\cdots,[a_r,b_r])$ is a snake with $a_1<\cdots<a_r$, $b_1<\cdots<b_r$ and $a_1=1-r$, then setting
$$\lambda=(b_r\ge b_{r-1}+1\ge \cdots\ge b_1+r-1),\ \ \mu=(a_r\ge a_{r-1}+1\ge\cdots\ge  a_2+r-2\ge 0)$$ defines a skew diagram  such that $\bor(\lambda/\mu)=\bos$.
Every snake is of this form up to translating all segments by the same integer; moreover this operation commutes with $\bos\mapsto\bos^\dagger$ and so there is no loss of generality in restricting our attention to such snakes. On the side of diagrams this corresponds to skew diagrams  with $\mu_i<\lambda_i$ for all $1\le i\le \ell$ and  $\mu_\ell=0$. The following lemma is stated in an equivalent form, without proof, in \cite[\S3.2]{LM14} and we include a proof for the reader's convenience.
\begin{lem}\label{skewdual} For every skew Young diagram $\lambda/\mu$ we have
$\bor(\lambda/\mu)^\dagger=\boc(\lambda/\mu).$
In particular, the Zelevinsky dual of a snake is a snake.
\end{lem}
\begin{pf} The second statement is immediate from the discussion preceding the lemma. For the first,  we proceed by induction on the number $\sum_i(\lambda_i-\mu_i)$. We may assume that $\mu_1<\lambda_1$, since
removing an empty first row translates both $\bor(\lambda/\mu)$ and $\boc(\lambda/\mu)$ by $1$. Order the segments of $\bor(\lambda/\mu)$ as in
Section \ref{ms}, and let $m$ be maximal such that $\lambda_1=\cdots=\lambda_m$. Since $\mu_s\le\mu_1<\lambda_1$ for $s\le m$, an inspection of
$\bor(\lambda/\mu)$ shows that its first $m$ terms are $[\mu_s-s+1,\lambda_1-s+1]$, $1\le s\le m$, and that the right end point of every other
segment is less than $\lambda_1-m$; hence \eqref{lead} gives
$$\bor(\lambda/\mu)^\dagger=([\lambda_1-m,\lambda_1])+\bor(\tilde\lambda/\mu)^\dagger,\ \ \tilde\lambda=(\lambda_1-1\ge\cdots\ge\lambda_1-1\ge\lambda_{m+1}\ge\cdots\ge\lambda_\ell).$$
An inspection shows that $\boc(\lambda/\mu)=([\lambda_1-m,\lambda_1])+\boc(\tilde\lambda/\mu)$, and the first statement of the lemma is immediate
from the induction hypothesis. 
\end{pf}
\vskip12pt
\noindent If $\lambda/\mu$ has no empty rows and no empty columns, Theorem \ref{altdet} gives
$${\rm sch}(\bor(\lambda/\mu))=\det(e_{\mu_i-i+1,\lambda_j-j+1})_{1\le i,j\le\ell},\qquad {\rm sch}(\boc(\lambda/\mu))=\det(e_{i-\lambda'_i,\,j-\mu'_j})_{1\le i,j\le\lambda_1}.$$
Together with Lemma \ref{skewdual}, Theorem \ref{zdualsoc}(ii) and the N\"agelsbach--Kostka identity \cite[(I.5.5)]{Mac95} we get:
\begin{cor}\label{skewphi} Let $\lambda/\mu$ be a skew Young diagram with $n\ge n(\bor(\lambda/\mu))$. Then
\begin{gather*}\Phi({\rm sch}(\bor(\lambda/\mu)))={\rm sch}(\boc(\lambda/\mu)),\ \ \kappa({\rm sch}(\boc(\lambda/\mu)))=s_{\lambda/\mu},\qquad \kappa({\rm sch}(\bor(\lambda/\mu)))=s_{\lambda'/\mu'}.\end{gather*}\hfill\qedsymbol
\end{cor}
\begin{rem}\label{latticerem} The elements $e_{a,b}$, $h_{a,b}$ are the elementary and complete functions of Macdonald's ninth
variation of Schur functions \cite{Mac92} (under $h_{rs}\leftrightarrow h_{s,s+r}$, with $e_{a,b}=0$ for $b-a>n+1$). For a snake, the
determinantal formula is thus a ninth variation Jacobi--Trudi formula, and Corollary \ref{skewphi} gives the corresponding
N\"agelsbach--Kostka identity; interpreting $e_{a,b}$ as a generating function of lattice paths, the Lindstr\"om--Gessel--Viennot
lemma \cite{Lin73,GV85} recovers the tableaux and path descriptions of \cite{BR90,MY12}, as in \cite{NN06}.
\end{rem}

\subsection{Alternating snakes and symmetric functions}\label{altclassical}
We now consider the symmetric functions $\kappa({\rm sch}(\bos))$ attached to alternating snakes, and what is possibly new on the
classical side. Throughout, $\bos=([i_1,j_1],\cdots,[i_r,j_r])$ is a stable alternating snake and $A(\bos)$ is the matrix of Section
\ref{snakesdef}. 
Since $\Phi_n(e_{a,b})=[V(\bomega_{a,b})]$ we may and do regard $A(\bos)$ as a matrix with entries in
${\rm Sym}(\mathbb I_n)$, and then we have that ${\rm sch}(\bos)=\det A(\bos)$.

\subsubsection{} Let $\bos=([i_1,j_1],\cdots,[i_r,j_r])$ be a stable alternating snake and assume $n\ge n(\bos)$. In
this range two segments of $\bos$ are connected if and only if they overlap, so the matrix $A(\bos)$ does not depend on $n$, and
$$\kappa({\rm sch}(\bos))=\det\kappa(A(\bos)).$$
Here $\kappa(A(\bos))_{p,\ell}$ is either $e_{j_\ell-i_p}$ or $0$, the zeros being those of $A(\bos)$. For an alternating snake which is not a snake,
$\kappa(A(\bos))$ has zeros in positions where $0\le j_\ell-i_p\le n+1$. For example,
if $\bos=([0,5],[2,6],[4,9],[3,4])$ then
$$\kappa(A(\bos))=\begin{pmatrix}e_5&e_6&e_9&0\\ e_3&e_4&e_7&0\\ e_1&e_2&e_5&1\\ e_2&e_3&e_6&e_1\end{pmatrix},$$
whereas the N\"agelsbach--Kostka matrix of the same end points would have $e_4$ and $e_2$ in place of the two zeros.

\subsubsection{} We say that an alternating snake $\bos$ is prime if $V(\bomega_\bos)$ is prime. By
\cite[Theorem 1]{BC25a}, $\bos=([i_1,j_1],\cdots,[i_r,j_r])$ is prime if and only if $([i_s,j_s],[i_{s+1},j_{s+1}])$ is connected for $1\le s\le r-1$
and $[i_{s-1},j_{s-1}]$ and $[i_{s+1},j_{s+1}]$ have no common end point for $2\le s\le r-1$. In general we can write (see \cite[Theorem 1]{BC25a})
$$V(\bomega_\bos)\cong V(\bomega_{\bos(p_0,p_1)})\otimes\cdots\otimes V(\bomega_{\bos(p_{k-1},p_k)}),\ \ 0=p_0<\cdots<p_k=r$$ where the $\bos(p_i,p_{i+1})$ are prime alternating snakes.

\begin{prop}\label{salt} Let $\bos$ be a stable alternating snake with $n\ge n(\bos)$, and let $\bos(p_i,p_{i+1})$, $0\le i\le k-1$, be its prime
factors, as in the preceding discussion.
\begin{enumerit}
\item[(i)] We have ${\rm sch}(\bos)={\rm sch}(\bos(p_0, p_1))\cdots{\rm sch}(\bos(p_{k-1}, p_k))$.
\item[(ii)] If $\bos$ is prime, let $\bos_p$, $1\le p\le r_1$, be the stable alternating snakes of \cite[Section 6]{BC25a}. Then
\begin{gather*}{\rm sch}(\bos)=\sum_{p=1}^{r_1}(-1)^{p+1}e_{i_p,j_1}\,{\rm sch}(\bos_p)\ \ \text{if }\ \ \epsilon_1=1,\ \
{\rm sch}(\bos)=\sum_{p=1}^{r_1}(-1)^{p+1}e_{i_1,j_p}\,{\rm sch}(\bos_p)\ \ \text{if $\epsilon_1=0$}.\end{gather*}
\item[(iii)] We have \begin{gather*}\kappa({\rm sch}(\bos))\in\sum_\lambda\mathbb Z_{\ge0}\,s_\lambda,\ \ \rho_{n+1}(\kappa({\rm sch}(\bos)))=\ch V(\bomega_\bos),\\ \kappa({\rm sch}(\bos^\dagger))=\omega(\kappa({\rm sch}(\bos)))=\det\omega(\kappa(A(\bos))),\ \ \rho_{n+1}(\kappa({\rm sch}(\bos^\dagger)))=\ch V(\bomega_{\bos^\dagger}).\end{gather*} 
\end{enumerit}
where  $\omega$ is applied to
$\kappa(A(\bos))$ entry-wise, i.e. each $e_k$ is replaced by $h_k$.
\end{prop}
\begin{pf}
    The first two statements are immediate from \cite[Theorem 3, Propositions 6.1--6.3]{BC25a}, Lemma \ref{kappa} and
\eqref{omegas}.  For (iii), the matrix $A(\bos)$ is essentially the same for all $n\ge n(\bos)$ and hence  $\kappa({\rm sch}(\bos))=\det\kappa(A(\bos))$ does not depend on $n$.  Taking $n\ge|\bos|-1$ we see that  the coefficients which appear in the Schur basis are multiplicities which gives non-negativity as needed.
\end{pf}

\begin{rem}\label{zeroexp} Writing each zero entry of $\kappa(A(\bos))$ as $e_{j_\ell-i_p}-e_{j_\ell-i_p}$ and expanding by
multilinearity, each $\kappa({\rm sch}(\bos))$ becomes a signed combination of products of elementary symmetric functions and
determinants $(e_{j-i})$ with increasing rows and columns, i.e. of skew Schur functions. For example, for
$\bos=([0,5],[2,6],[4,9],[3,4])$ as above,
$$\kappa({\rm sch}(\bos))=\kappa({\rm sch}(\boT_0))-e_4\,\kappa({\rm sch}(\boT_1))+e_2\,\kappa({\rm sch}(\boT_2)),$$
where $\boT_0=\{[0,4],[2,5],[3,6],[4,9]\}$, $\boT_1=\{[2,5],[3,6],[4,9]\}$ and $\boT_2=\{[0,5],[3,6],[4,9]\}$ are snakes, so that
each $\kappa({\rm sch}(\boT_k))$ is a skew Schur function (the N\"agelsbach--Kostka formula). This expression is not positive and does
not explain the Schur positivity of $\kappa({\rm sch}(\bos))$.
The determinantal formula for alternating snakes  resembles the situation in types
$C_n$ and $D_n$ studied in \cite{NN06,NN07}.
\end{rem}

\begin{rem}\label{bittmann} Zeros of a similar nature occur in the work of Bittmann \cite{Bit23} which was discussed in the introduction and part (iii) of the preceding proposition applies for these examples as well.
\end{rem}

\section{Simplicity of ${\rm hd}\, M(\bos)$ and $\soc M(\bos)$}\label{head}
We begin this section by translating some well-known results about $\mathscr F_n$ into the language of intervals. We also recall some results from \cite{KKOP}. We then establish the simplicity of the head and socle of $M(\bos)$, and conclude the section with further representation theoretic results which are needed to prove part (ii) of Theorem \ref{zdualsoc}.

\subsection{} Given an object $V$ of $\Fn$ we denote by $V^*$ the right dual of $V$, i.e., we have an inclusion $\mathbb C\hookrightarrow V^*\otimes V$. Given $\bos\in\mult(\In)$ we define $\bos^*\in\mult(\In)$ by giving its multiplicity function as follows:
$$\mu_{\bos^*}([d,c+n+1])=\mu_\bos([c,d]).$$
It was proved in \cite{CP95a} that $V(\bomega_\bos)^*\cong V(\bomega_{\bos^*})$.
\subsection{} A module $V$ is said to be $\ell$-highest weight if it is a quotient of $W(\bomega_{\bos})$ for some $\bos\in\mult(\In)$. In particular, an $\ell$-highest weight module has a unique irreducible quotient, i.e., it has a simple head. Recall that a module is real if its tensor square is irreducible.
\\\\
We recall special cases of a result from \cite{ch01} and translate them into the language of intervals.
\begin{prop} \label{chb}
Suppose that $([a_1,b_1], [a_2,b_2],\cdots, [a_r,b_r])$ is a sequence of intervals in $\mathbb I_n^+$ such that $b_1\ge\cdots\ge b_r$. \begin{enumerit} \item[(i)] The modules $$V(\bomega_{a_1,b_1})\otimes\cdots\otimes V(\bomega_{a_r,b_r}),\ \ V(\bomega(a_1,b_1))\otimes\cdots\otimes V(\bomega(a_r,b_r))$$ are $\ell$-highest weight.
\item[(ii)] Let $[a_1,b_1],[a_2,b_2]\in\mathbb I_n^+$. If the pair $([a_1,b_1],[a_2,b_2])$ is not connected (resp. not overlapping) then the module $$V(\bomega_{a_1,b_1})\otimes V(\bomega_{a_2,b_2}),  \ \ \ {\rm(resp.}\ \  V(\bomega(a_1,b_1))\otimes V(\bomega(a_2,b_2)) {\rm )}$$ is irreducible.
        If  the pair is connected (resp. overlapping), say, $a_1<a_2\leq b_1<b_2$  we have short exact sequences
        \begin{gather*}0\to V(\bomega_{a_1,b_2}\bomega_{a_2,b_1})\to V(\bomega_{a_2,b_2})\otimes V(\bomega_{a_1,b_1})\to V(\bomega_{a_1,b_1}\bomega_{a_2,b_2})\to 0,\\
        0\to V(\bomega')\to V(\bomega(a_2,b_2))\otimes V(\bomega(a_1,b_1))\to V(\bomega(a_1,b_2)\bomega(a_2,b_1))\to 0,
        \end{gather*}
        where $\bomega'=\bomega(a_1,a_2-1)\bomega_{a_2-1, a_2+1}\cdots\bomega_{b_1-1,b_1+1}\bomega(b_1+1,b_2)$.
        \item[(iii)] The module $V(\bomega(a,b))^*\otimes  V(\bomega(c,d))$ is irreducible if $n\ge d-\min\{a,c\}$ and either $b<d$, or $b=d$ and $a<c$.\hfill\qedsymbol
\end{enumerit}
\end{prop}
\subsection{} In \cite{KKOP}, the authors defined, for simple modules $M,N$, a non-negative integer $\lie d(M, N)$ depending on $n$, which is the order of the zero at $z=1$ of $d_{M,N}(z)d_{N,M}(z^{-1})$, where $d_{M,N}(z)$ denotes the denominator of the normalized $R$-matrix of $M$ and $N$. We summarize some important properties of $\lie d$.
\begin{prop}\label{kkop}
Let $L$, $M$ and $N$ be simple modules. Then:
\begin{enumerit}
   \item[(i)] $\lie d (M,N)= \lie d (N,M) = \lie d (M^*,N^*)$.
   \item[(ii)] If  $M\otimes N$ is irreducible then $$\lie d(L,M\otimes N) = \lie d(L,M) +\lie d(L, N).$$
   \item[(iii)] For any simple subquotient $S$ of $M\otimes N$ we have
   $$\lie d(L,S)\leq \lie d(L,M)+ \lie d(L, N).$$
   \item[(iv)]  Suppose that $M$ is real.  Then $M\otimes N$ is irreducible if and only if $\lie d (M,N)=0$.
   \item[(v)] Suppose that $M_1,\cdots, M_k$ are real modules such that $\lie d (M_i^*, M_j)=0$ if $i>j$. Then $M_1\otimes\cdots\otimes M_k$ has simple socle and
   $$\soc (M_1\otimes \cdots \otimes M_k)\cong {\rm hd}(M_k\otimes \cdots\otimes M_1).$$

\end{enumerit}\hfill\qedsymbol
\end{prop}
We record  the following trivial consequences of the proposition.
\begin{cor} \label{tensormix}
\begin{enumerit}
    \item[(i)] Suppose that $M_1,\cdots, M_k$ are real simple modules such that $\lie d(M_i,M_j)=0$ for $1\le i,j\le k$. Then $M_1\otimes \cdots\otimes M_k$ is real and simple.
    \item[(ii)] Let $N, M_1,\cdots, M_k$ be real modules such that
    $$\lie d (M_i^*,M_j)=0, \ \ i>j, \ \ \lie d(M_i^*,N)=0, \ \ i\ge 2, \ \ \ \lie d(M_1,N)=0.$$
   Then $M_1\otimes N\otimes M_2\otimes\cdots\otimes M_k$ has simple socle.
\end{enumerit}
\end{cor}

\subsection{Proof of Theorem \ref{zdualsoc}(i) and simplicity of $\soc M(\bos)$} \label{pfi} Proposition \ref{chb}(ii) shows that the modules $V(\bomega(a,b))$ are real for all $[a,b]\in\In$.\\\\
By our conventions a multisegment $\bos=\{[a_1,b_1],\cdots, [a_r, b_r]\}$ is ordered so that $b_1\ge\cdots\ge b_r$. It is then immediate from Proposition \ref{chb}(i) that $M(\bos)$ is $\ell$-highest weight and hence has simple head $V(\bomega(\bos))$. The definition of $\bos^\sharp$ and \eqref{sig} show that $\bomega(\bos)=\bomega(\bos^\sharp)$. Since $M(\bos^\sharp)$ is $\ell$-highest weight by Proposition \ref{chb}(i), $V(\bomega(\bos))$ is also a quotient of $M(\bos^\sharp)$, and hence part (i) of the theorem follows if we prove that $M(\bos^\sharp)$ is irreducible.
\\\\
By construction (see Section \ref{std}) the multisegment $\bos^\sharp$ consists of intervals which are pairwise not connected. Hence Proposition \ref{chb}(ii), Proposition \ref{kkop}(iv) and Corollary \ref{tensormix}(i) prove that $M(\bos^\sharp)$ is irreducible, as needed. \\\\
 To prove that $M(\bos)$ has simple socle when $n\ge n(\bos)$, we fix
an enumeration  $c_1>\cdots>c_k$ of the distinct elements of the multiset $\{b_1,\cdots, b_r\}$ and set
$$M_s= \bigotimes_{\{p: b_p=c_s\}} V(\bomega(a_p, c_s)),\ \ 1\le s\le k.$$ Proposition \ref{chb}(ii) and Proposition \ref{kkop}(iv) show that $\lie d(V(\bomega(a_p, c_s)), V(\bomega(a_\ell, c_s)))=0$ for all $p,\ell$ with $b_p=b_\ell=c_s$. It follows from Corollary \ref{tensormix}(i) that $M_s$ is real and simple. \\\\ Proposition \ref{chb}(iii) and Proposition \ref{kkop}(iv) give $$\lie d( V(\bomega(a_s, b_s))^*, V(\bomega(a_p,b_p)))=0\ \ \text{if}\ \ s>p\ \ {\rm and}\ \ b_s<b_p\ \ ({\rm i.e.}\ \ b_s\ne b_p). $$ Hence a repeated application of  Proposition \ref{kkop}(i),(ii) gives
$$\lie d(M_s^*, M_p)=0,\ \ 1\le p<s\le k. $$ Proposition \ref{kkop}(v) now shows that $M(\bos)$ has simple socle.

\subsection{} We need the following result in the following sections. Here we are under the blanket assumption that none of the modules involved in the statement are trivial. \begin{prop}\label{unmixedpair} \begin{enumerit}
    \item[(i)] For $n\ge b-c$  we have $$\lie d(V(\bomega_{a,b})^*, V(\bomega(c,d)))=0\ {\rm if}\  d<b,\ \  \lie d( V(\bomega(c,d))^*, V(\bomega_{a,b}))=0 \ {\rm if}\  d\le a.$$
    \item[(ii)] If  $n\geq b-\min\{a,c\}$ then
    $$\lie d (V(\bomega_{a,b})^*,V(\bomega(c,b)))=1,\ \ {\rm and}\ \ \lie d (V(\bomega_{a,b})^*,V(\bomega_{c,b}))=1\ \ {\rm if}\ \ a\le c.$$
    \item[(iii)] If $d<a$, then $V(\bomega_{a,b})\otimes V(\bomega(c,d))\cong V(\bomega(c,d))\otimes V(\bomega_{a,b})$
    is irreducible.
\end{enumerit}
\end{prop}
    \begin{proof}  We use freely the fact that $V(\bomega_{a,b})^* \cong V(\bomega_{b,a+n+1})$. For the first equality (resp. the second equality) in part (i)
 we note that our conditions imply that $[b,a+n+1]$ and $[j,j+1]$ are not connected for all $c\le j<d$ (resp. $[j+1,j+n+1]$ and $[a,b]$  are not connected for all $c\le j<d$). Hence Propositions \ref{chb}(ii) and \ref{kkop}(iv) give $$\lie d(V(\bomega_{a,b})^*, V(\bomega_{j,j+1}))=0=\lie d(V(\bomega_{j,j+1})^*, V(\bomega_{a,b})),\ \ c\le j<d.$$
Part (i)  is now  immediate from a repeated application of  Proposition \ref{kkop}(iii).\\\\
To prove (ii), we note that $([b-1,b],[b,a+n+1])$  is connected and hence $([c,c+1],\cdots, [b-1,b],[b,a+n+1])$ is a prime snake and the first equality is a consequence of  \cite[Lemma 6.2.1]{Naoi24}.
The second part is immediate from Proposition \ref{chb}(ii) since our conditions imply that
the intervals $([b,a+n+1],[c,b])$ are connected.\\\\
For (iii), note that if $c\le e<d$ then $e+1<a$, so that the pair $([e,e+1],[a,b])$ is not connected,
and Propositions \ref{chb}(ii) and \ref{kkop}(iv) give $\lie d(V(\bomega_{a,b}),V(\bomega_{e,e+1}))=0$.
Since $V(\bomega(c,d))$ is a quotient of $V(\bomega_{d-1,d})\otimes\cdots\otimes V(\bomega_{c,c+1})$,
Proposition \ref{kkop}(iii) gives $\lie d(V(\bomega_{a,b}),V(\bomega(c,d)))=0$, and (iii) follows from
Proposition \ref{kkop}(iv) since $V(\bomega_{a,b})$ is real.
  \end{proof}
\section{The isomorphism {$\soc (M(\bos))\cong V(\bomega_{\bos^\dagger})$}}\label{socMs}
Throughout this section we assume that $\bos=\{[c_1,d_1],\cdots, [c_r,d_r]\}\in\mult(\In)$ with $n\ge n(\bos)$ and we let $([a_1,b_1],\cdots, [a_m,b_1-m+1])$ and $\bos[1]$ be  as
in Section \ref{zdual}. Clearly the multiplicity functions $\mu_\bos$ and $\mu_{\bos[1]}$ are related by
$$\mu_{\bos[1]}([c,d])=\mu_\bos([c,d])-\sum_{p=1}^m\delta_{[c,d],[a_p,b_1-p+1]}+\sum_{p=1}^m\delta_{[c,d],[a_p,b_1-p]}.$$

\subsection{} The first  step is to define
multisegments $\bos_p$ and $\bos'_p$  for $1\le p\le m$ as follows: writing $\mu_{\bos_p}$ and $\mu_{\bos'_p}$ for their  multiplicity functions we set
\begin{gather*} \mu_{\bos_p}([c,d])=\begin{cases} 0,\ \ {\rm if}\ \ d<  b_1-p\ \ {\rm or}\ \ d=b_1-p\ \ {\rm and }\ \ c<a_p,\\
\mu_\bos([c,b_1-s])-\delta_{c,a_{s+1}}+\delta_{c,a_{s}},\ \ d=b_1-s,\ \ 0\le s<p,\\
\mu_\bos([c, b_1-p])+\delta_{c, a_p},\ \ d=b_1- p\ \ {\rm and}\ \ c\ge a_p,
\end{cases}
\\
\mu_{\bos'_p}([c,d])=\begin{cases} 0,\ \ {\rm if} \ d>b_1-p\ \ {\rm or}\ \ d=b_1-p\ \ {\rm and}\ \ c\ge a_p,\\ \mu_\bos([c,d]) - \delta_{[a_{p+1},b_1-p],[c,d]},\ \ d< b_1-p,\ \ {\rm or} \ d=b_1-p\ \ {\rm and}\ \  c<a_p.\\
   \end{cases}
\end{gather*}
Here we understand that $a_0=\infty$ and  $\delta_{[a_{m+1},b_1-m],[c,d]}=0$. \\\\
 We record the following for later use.

\begin{lem}\label{mpcount}
We have $$\displaystyle\sum_{c\in \mathbb Z}\mu_{\bos_p}([c,b_1])=\sum_{c\ge b_1-p}\mu_{\bos_p^\dagger}([c,b_1]).$$
\end{lem}
\begin{pf} It is clear from \eqref{lead} that
$$N:=\displaystyle\sum_{c\in \mathbb Z}\mu_{\bos_p}([c,b_1])=\sum_{c\in\mathbb Z}\mu_{\bos_p^\dagger}([c,b_1]).$$
Let $\bos_p[k]$, $0\le k\le N$, be such that
\begin{gather*}
(\bos_p[k])^\dagger=([b_1-m_{k+1}, b_1])+ (\bos_p[k+1])^\dagger,\ \ \bos_p[0]=\bos_p,\\
\bos_p^\dagger=([b_1-m_1, b_1],\cdots, [b_1-m_N, b_1])+ (\bos_p[N])^\dagger.
\end{gather*}
The lemma is equivalent to the statement that $m_k\le p$ for $1\le k\le N$.
\\\\
For $1\le k\le N$ let $\ell_{1,k}>\ell_{2,k}>\cdots>\ell_{m_k,k}$ be the chain defining $\bos_p[k]$
from $\bos_p[k-1]$ as in \eqref{lead}, so that
\begin{gather*}
\ell_{1,k}=\max\{c: \mu_{\bos_p[k-1]}([c,b_1])>0\}, \\
\ell_{j,k}=\max\{c:\ \mu_{\bos_p[k-1]}([c,b_1-j+1])>0,\ c<\ell_{j-1,k}\},\qquad 2\le j\le m_k.
\end{gather*}
We first show that $m_{k+1}\le m_k$. Since the segments of $\bos_p[k]$ with right end point $b_1$
are segments of $\bos_p[k-1]$, we have $\ell_{1,k+1}\le\ell_{1,k}$. Suppose that $\ell_{j-1,k+1}\le\ell_{j-1,k}$ for some $2\le j\le m_{k+1}$. Then
$\ell_{j,k+1}<\ell_{j-1,k+1}\le\ell_{j-1,k}$, so $[\ell_{j,k+1},b_1-j+1]$ is not the segment
$[\ell_{j-1,k},b_1-j+1]$ obtained by the recursive process, and hence is a segment of $\bos_p[k-1]$. If $j\le m_k$ the definition of $\ell_{j,k}$ gives
$\ell_{j,k+1}\le\ell_{j,k}$. If $j=m_k+1$, the chain $\ell_{1,k}>\cdots>\ell_{m_k,k}$ in
$\bos_p[k-1]$ could be extended by $\ell_{j,k+1}$, contradicting the maximality of $m_k$. Hence
$m_{k+1}\le m_k$.
\\\\
It remains to prove that $m_1\le p$. The same argument, comparing $\bos_p$ with $\bos$ and using
that the only segment of $\bos_p$ with right end point $b_1-j+1$ which is not a segment of $\bos$ is
$[a_{j-1},b_1-j+1]$, shows that $\ell_{j,1}\le a_j$ for $1\le j\le\min\{m_1,p\}$. If $m_1\ge p+1$,
then $[\ell_{p+1,1},b_1-p]$ is a segment of $\bos_p$, so that $\ell_{p+1,1}\ge a_p$ by the definition
of $\bos_p$, while $\ell_{p+1,1}<\ell_{p,1}\le a_p$. Hence $m_1\le p$.
\end{pf}

\subsection{}\label{mp} Let $\bom_p$ be the multisegment whose multiplicity function is given by $$\mu_{\bom_p}([c,d])=\begin{cases}\mu_\bos([c,b_1-p+1])-\delta_{c,a_p}, & d=b_1-p+1,\ c<a_{p-1},\\ \mu_\bos([c,b_1-p]), & d=b_1-p,\ c\ge a_p,\\ 0,&\text{otherwise.}\end{cases}$$

It is easy to check that for $1\le p\le m$
\begin{equation}\label{sprel}\bos_p=\bos_{p-1} + \bom_p + ([a_p,b_1-p]),\ \ \ \bos'_{p-1}=\bom_p + ([a_{p+1},b_1-p]) +\bos'_p,
\end{equation}
where $\bos_0=\emptyset$ and the middle term is absent in $\bos_{m-1}'$.
We have,

\begin{lem}\label{mplem}
\begin{gather*}M(\bos_p)\cong M(\bos_{p-1})\otimes M(\bom_p)\otimes V(\bomega(a_p,b_1-p)),\\
M(\bos_{p-1}')\cong M(\bom_p)\otimes V(\bomega(a_{p+1},b_1-p))^{1-\delta_{p,m}}\otimes M(\bos'_p).
\end{gather*}
\end{lem}
\begin{pf}Assume that $\bos_{p-1}$ and $\bom_p$ are written in the lexicographic order of Section
\ref{ms}. An inspection shows that the expression in \eqref{sprel} respects this order, and hence
the first isomorphism is immediate by definition. Noting the  choice of $a_{p+1}$, an identical proof gives the second isomorphism.
\end{pf}

\subsection{}\label{induction}  We prove that $\soc (M(\bos))\cong V(\bomega_{\bos^\dagger})$ by induction on $|\bos|=\sum_{s=1}^r(d_s-c_s)$. Clearly induction begins at $|\bos|=1$. Assume now that the result is proved for all multisegments $\bos'$ with $|\bos'|<|\bos|$ and
$n\ge n(\bos')$.
Define modules
\begin{gather*} U_0= M(\bos),\ \ U_p= M(\bos_p)\otimes V(\bomega_{b_1-p-1,b_1}\bomega(a_{p+1}, b_1-p-1))\otimes M(\bos_{p}'),\ \ 1\le p\le m-1,\\ W_p= M(\bos_p)\otimes V(\bomega_{b_1-p,b_1})\otimes V(\bomega(a_{p+1},b_1-p))^{1-\delta_{p,m}}\otimes  M(\bos_p'), \ \ 1\le p\le m.
\end{gather*}
We shall prove that
\begin{gather}\label{uw}
U_p\hookrightarrow W_p\ \ 1\le p\le m-1,\ \  U_{p-1}\hookrightarrow W_{p},\ \ \  \soc W_p\cong V(\bomega_{\bos^\dagger}),\ \ \ 1\le p\le m.
\end{gather}

 \subsection{Proof of the inclusions
$U_{p-1}\hookrightarrow W_p\hookleftarrow U_p$} We note first that for all $1\le p\le m-1$ the multisegment $([a_{p+1},a_{p+1}+1],\cdots, [b_1-p-1,b_1-p],[b_1-p,b_1])$ is a prime snake. Moreover  the module $V(\bomega_{b_1-p,b_1})\otimes V(\bomega(a_{p+1}, b_1-p))$ is clearly a quotient of the corresponding standard module and hence is $\ell$-highest weight.   The following identity in $\cal K(\Fn)$ was proved in
 \cite[Proposition 9]{BC25a},
 \begin{gather*}
[V(\bomega_{b_1-p,b_1})][V(\bomega(a_{p+1}, b_1-p))]=    [V(\bomega_{b_1-p,b_1}\bomega(a_{p+1}, b_1-p))]+[V(\bomega_{b_1-p-1,b_1}\bomega(a_{p+1},b_1-p-1))].
 \end{gather*}
 It now follows that   for $1\leq p\leq m-1$, we have
    $$V(\bomega_{b_1-p-1,b_1}\bomega(a_{p+1},b_1-p-1))\cong \soc(V(\bomega_{b_1-p,b_1})\otimes V(\bomega(a_{p+1},b_1-p))).$$ Tensoring on the left with $M(\bos_p)$ and on the right with $M(\bos'_p)$ gives $U_p\hookrightarrow W_p$.\\\\
    We claim  that the inclusion $U_{p-1}\hookrightarrow W_p$ can be deduced from the following
\begin{equation} \label{mpcom}V(\bomega_{b_1-p,b_1}\bomega(a_p,b_1-p))\otimes M(\bom_p)\cong M(\bom_p)\otimes  V(\bomega_{b_1-p,b_1}\bomega(a_p,b_1-p)).\end{equation}
For the claim, we set $\bomega_p= \bomega_{b_1-p,b_1}\bomega(a_p, b_1-p)$ and use  the definition of $U_{p-1}$, Lemma \ref{mplem} and the isomorphism in \eqref{mpcom} to get
\begin{gather}\label{up} U_{p-1}\cong M(\bos_{p-1})\otimes M(\bom_p)\otimes V(\bomega_p)\otimes V(\bomega(a_{p+1},b_1-p))^{1-\delta_{p,m}}\otimes M(\bos'_p).\end{gather}
(For $p=1$ this is equivalent to saying $U_0=M(\bos)$ since $\bos_0=\emptyset$ and $V(\bomega_1)=V(\bomega(a_1,b_1))$.) \\\\
The module $V(\bomega_{b_1-p,b_1})\otimes V(\bomega(a_p,b_1-p))$ is clearly a quotient of the corresponding standard module and hence is $\ell$-highest weight with simple head $V(\bomega_p)$. Since $V(\bomega_{b_1-p,b_1})$ is real, it follows from \cite[Theorem 3.12]{KKKO15} that
$V(\bomega_p)\hookrightarrow V(\bomega(a_p,b_1-p))\otimes V(\bomega_{b_1-p,b_1})$, and so  the isomorphism in \eqref{up} shows that $U_{p-1}$ is included in
\begin{gather*}
 M(\bos_{p-1})\otimes M(\bom_p)\otimes V(\bomega(a_p,b_1-p))\otimes V(\bomega_{b_1-p,b_1})\otimes V(\bomega(a_{p+1},b_1-p))^{1-\delta_{p,m}}\otimes M(\bos'_p).
\end{gather*}
An application of Lemma \ref{mplem} now shows that the preceding module is isomorphic to $W_p$ as needed.
\\\\
We prove the isomorphism \eqref{mpcom}.
  The multisegment $\bom_p$ consists of intervals of the form $[c,b_1-p+1]$ with $c<a_{p-1}$, and hence $c\leq a_p$ by the choice of $a_p$, and $[d,b_1-p]$ with $a_p\leq d$. Since $V(\bomega_p)$ is a snake module, it is real, and an application of Proposition \ref{kkop}(iv) shows that it suffices to prove
\begin{gather}
\label{opcom1} \lie d (V(\bomega_p), V(\bomega(d,b_1-p)))=0, \ \  a_p\leq d,\ \
\lie d (V(\bomega_p), V(\bomega(c,b_1-p+1)))=0, \ \ c\leq a_p.
\end{gather}
 For the first equality    we note that  $\bomega_p$ is given by the multisegment $([a_p,a_p+1],\cdots,[b_1-p-1,b_1-p], [b_1-p,b_1])$ which is  a snake. Moreover $\bomega(d,b_1-p)$ divides $\bomega_p$ if $a_p\le d$ and hence the result follows from \cite[Proposition 4.1.3(ii)]{Naoi24}. \\\\
For the second equality in \eqref{opcom1}, a further application of Proposition \ref{kkop}(iii) shows that it suffices to prove
\begin{gather}\label{opcom2} \lie d (V(\bomega(a_p,b_1-p)), V(\bomega(c,b_1-p+1))) =0= \lie d (V(\bomega_{b_1-p,b_1}), V(\bomega(c,b_1-p+1))).\end{gather}
Since $c\leq a_p\leq b_1-p<b_1-p+1$ the pair $([a_p,b_1-p],[c,b_1-p+1])$ is not overlapping and the first equality is immediate from Propositions \ref{chb}(ii) and \ref{kkop}(iv).\\\\
 The second equality is also immediate if  $c=b_1-p$ since the pair $([b_1-p,b_1-p+1], [b_1-p,b_1])$ is not connected and  Propositions \ref{chb}(ii) and \ref{kkop}(iv) apply. If $c<b_1-p$, we have
\begin{gather*}
    \lie d (V(\bomega(c,b_1-p+1)), V(\bomega_{b_1-p,b_1}))\leq \\ \lie d (V(\bomega(c,b_1-p-1)), V(\bomega_{b_1-p,b_1})) + \lie d (V(\bomega(b_1-p-1,b_1-p+1)), V(\bomega_{b_1-p,b_1})).
\end{gather*}
Since $b_1-p-1<b_1-p$, Proposition \ref{unmixedpair}(iii) and Proposition \ref{kkop}(iv) show that the first term on the right hand side of the inequality is zero. \\\\
For the second term we note that the multisegment $([b_1-p-1, b_1-p], [b_1-p,b_1-p+1], [b_1-p,b_1])$ is a non-prime alternating snake.
By  \cite[Theorem 1]{BC25a} the associated prime decomposition of the multisegment is  $([b_1-p,b_1])+ ([b_1-p-1,b_1-p],[b_1-p,b_1-p+1])$, or equivalently the module $V(\bomega(b_1-p-1, b_1-p+1))\otimes V(\bomega_{b_1-p,b_1})$ is irreducible; hence the second term on the right hand side is also zero by Proposition \ref{kkop}(iv), and the proof of \eqref{opcom2} and hence of the inclusions is complete.
 \subsection{Normal sequences} To establish the simplicity of the socle of the modules $W_p$ we need some additional definitions and results from \cite{KKOP}.
 \subsubsection{}  Say that a sequence $(M_1,\cdots, M_k)$ of real simple modules is normal if the composition of the $R$-matrices
$$\bor_{M_{k-1},M_k}\circ \cdots \circ (\bor_{M_1,M_k}\circ\cdots\circ\bor_{M_1,M_2}): M_1\otimes \cdots \otimes M_k \to M_k \otimes \cdots\otimes M_1$$
is non-zero. Say that a pair $(M,N)$ of real simple modules is unmixed if $\lie d(M^*,N)=0$, and that the sequence $(M_1,\cdots,M_k)$ is unmixed if $(M_i,M_j)$ is unmixed for all $1\le i<j\le k$.
\begin{prop}\cite[Lemma 4.15]{KKOP}\label{simplehead}
    If $(M_1,\cdots, M_k)$ is a normal sequence then $M_1\otimes \cdots\otimes M_k$ has simple head which is isomorphic to the socle of $M_k\otimes\cdots\otimes M_1$.
\end{prop}

 \subsubsection{} Given a module $M$ and $k\ge 0$, let $\mathscr D^0 M=M$ and let $\mathscr D^kM$  be the $k$-th right dual of $M$. Thus \begin{equation}\label{kdual} \mathscr D^{2k-1}V(\bomega_{a,b})\cong V(\bomega_{b+(k-1)(n+1), a+k(n+1)}),\ \ \ \mathscr D^{2k} V(\bomega_{a,b})\cong V(\bomega_{a+k(n+1), b+k(n+1)}).\end{equation}
If $M,N$ are simple modules, then it is known that $$\lie d(M,\mathscr D^kN)=0=\lie d(\mathscr D^kM,N)\ \ {\rm if}\ \ k\gg0,$$
and hence we have the following well-defined integers:
 \begin{gather*}
     \Lambda(M,N) = \sum_{k\in \mathbb Z_{\geq 0}}(-1)^k\lie d (M,\mathscr D^k N) - \sum_{k\in \mathbb Z_{> 0}}(-1)^{k}\lie d (\mathscr D^kM, N),\\
\Lambda^\infty(M,N) = \sum_{k\in \mathbb Z_{\ge 0}}(-1)^k\lie d (M,\mathscr D^k N)+ \sum_{k\in \mathbb Z_{> 0}}(-1)^k\lie d (\mathscr D^k M, N).
    \end{gather*}
    In particular, we have
   \begin{equation}\label{ldiff} \Lambda(M,N)-\Lambda^\infty(M,N)=2\sum_{k\in \mathbb Z_{> 0}}(-1)^{k+1}\lie d (\mathscr D^kM, N).\end{equation}
    It follows from Proposition \ref{kkop}(i)
    that
    $$\Lambda(M,N)=\Lambda(N, M^*),\ \ \Lambda^\infty(N,M)=\Lambda^{\infty}(M,N).$$
Finally, by \cite[Proposition 3.11]{KKOP}, for any simple subquotient $S$ of $M\otimes N$ we have \begin{equation}\label{ladd}\Lambda^{\infty}(L, S) = \Lambda^{\infty}(L, M) + \Lambda^{\infty}(L, N). \end{equation}  \subsubsection{Normality and simple heads} The following was proved in \cite[Lemma 5.3]{KKOP2} and \cite[Lemma 4.16]{KKOP}.
\begin{prop}\label{kkop1} Let $(M_1,\cdots,M_k)$ be a sequence of real simple modules.

\begin{enumerit}
\item[(i)] If $\lie d(M_i^*, M_j)=0$ for all $i<j$ then $(M_1,\cdots,M_k)$ is a normal sequence.
\item[(ii)] Suppose that $(M_2,\cdots,M_k)$ is normal. Then $(M_1,\cdots,M_k)$ is normal if and only if
$$\Lambda(M_1,{\rm hd}(M_2\otimes\cdots\otimes M_k))=\sum_{j=2}^k\Lambda(M_1,M_j).$$
\end{enumerit}
\end{prop}
\subsubsection{}  The following simplification will be helpful.\begin{prop}\label{Linfty}

\begin{enumerit}
\item[(i)] Let $\bos=([a_1,b_1],\cdots,[a_r,b_r])$ and $\bos'=([c_1,d_1],\cdots,[c_p,d_p])$ be multisegments such that
$n\geq \max\{b_i,d_j\}-\min\{a_i,c_j\}$. Then,
$$\Lambda(V(\bomega_{\bos}),V(\bomega_{\bos'})) = \Lambda^{\infty}(V(\bomega_{\bos}),V(\bomega_{\bos'})) +
2\lie d(V(\bomega_{\bos})^*,V(\bomega_{\bos'})).$$
\item[(ii)] Suppose that $\boT_1,\cdots, \boT_k\in\mult(\In)$ are such that $V(\bomega_{\boT_s})$ is real
for all $1\le s\le k$ and $n\ge n(\boT_1+\cdots+\boT_k)$. Suppose also that
$(V(\bomega_{\boT_2}),\cdots, V(\bomega_{\boT_k}))$ is a normal sequence and that
${\rm hd}(V(\bomega_{\boT_2})\otimes\cdots\otimes V(\bomega_{\boT_k}))\cong V(\bomega_{\boT_0})$. Then
$(V(\bomega_{\boT_1}),\cdots, V(\bomega_{\boT_k}))$ is normal if and only if
$$\lie d(V(\bomega_{\boT_1^*}), V(\bomega_{\boT_0}))=\sum_{j=2}^k\lie d(V(\bomega_{\boT_1^*}), V(\bomega_{\boT_j})).$$
\end{enumerit}
\end{prop}
\begin{pf}
The restriction on $n$ along with \eqref{kdual} and Propositions \ref{chb}(ii) and \ref{kkop}(iv) gives
$$\lie d(\mathscr D^kV(\bomega_{a_s,b_s}), V(\bomega_{c_t,d_t}))=0,\ \ k\ge 2,\ \ 1\le s\le r,\ \ 1\le t\le p.$$
An application of Proposition \ref{kkop}(iii) gives $\lie d(\mathscr D^k V(\bomega_\bos), V(\bomega_{\bos'}))=0$,
$k\geq 2$, and part (i) now follows from \eqref{ldiff}.
\\\\
Since $V(\bomega_{\boT_0})$ is a composition factor of $V(\bomega_{\boT_2})\otimes\cdots\otimes
V(\bomega_{\boT_k})$, it follows from the discussion in Section \ref{bou} that $V(\bomega_{\boT_0})$ is an object of
$\mathscr F_n(\bou_c)$, where $c=\min\{a:[a,b]\in\boT_1+\cdots+\boT_k\}$, and hence
$n\ge n(\boT_1+\boT_j)$ for $0\le j\le k$. A repeated application of
\eqref{ladd} gives
$\Lambda^\infty(V(\bomega_{\boT_1}),V(\bomega_{\boT_0}))=\sum_{j=2}^k\Lambda^\infty(V(\bomega_{\boT_1}),V(\bomega_{\boT_j}))$,
and hence part (i) gives
$$\Lambda(V(\bomega_{\boT_1}),V(\bomega_{\boT_0}))-\sum_{j=2}^k\Lambda(V(\bomega_{\boT_1}),V(\bomega_{\boT_j}))
=2\Big(\lie d(V(\bomega_{\boT_1^*}),V(\bomega_{\boT_0}))-\sum_{j=2}^k\lie d(V(\bomega_{\boT_1^*}),
V(\bomega_{\boT_j}))\Big).$$
Part (ii) is now immediate from Proposition \ref{kkop1}(ii).
\end{pf}
\subsection{A normal sequence for $W_p$} For  $1\leq p\leq m$ write
\begin{gather*}V_{p,k} = \bigotimes_{c\in \mathbb Z} V(\bomega(c,k))^{\otimes\mu_{\bos_p}([c,k])}, \ \ b_1-p\leq k\leq b_1,\\
V_{p,k}'= V(\bomega(a_{p+1},b_1-p))^{\otimes\delta_{b_1-p,k}(1-\delta_{p,m})}\otimes \bigotimes_{c\in \mathbb Z} V(\bomega(c,k))^{\otimes\mu_{\bos_p'}([c,k])}, \ \ d_r\le k\leq b_1-p.
\end{gather*}
Since the factors of $V_{p,k}$ (resp. $V'_{p,k}$) correspond to intervals with the same right end point $k$, Propositions \ref{chb}(ii), \ref{kkop}(iv) and Corollary \ref{tensormix}(i) show that  $V_{p,k}$ and $V_{p,k}'$ are real and simple modules.
It is easy to check that
$$M(\bos_p)\cong V_{p,b_1}\otimes \cdots \otimes V_{p,b_1-p} \ \ {\rm and} \ \ V(\bomega(a_{p+1},b_1-p))^{1-\delta_{p,m}}\otimes M(\bos_p') \cong V_{p,b_1-p}'\otimes \cdots \otimes V_{p,d_r}',$$
and hence
$$W_p \cong  V_{p,b_1}\otimes \cdots \otimes V_{p,b_1-p} \otimes V(\bomega_{b_1-p,b_1})\otimes V_{p,b_1-p}'\otimes \cdots \otimes V_{p,d_r}'.$$
Consider the sequence of real simple modules
$$\bs_p = (V_{p,d_r}', \cdots, V_{p,b_1-p}', V(\bomega_{b_1-p,b_1}), V_{p,b_1-p},\cdots, V_{p,b_1}).$$
We prove the following.
\begin{prop}\label{Spnormal}
    For $1\leq p\leq m$ the sequence $\bs_p$ is normal. In particular $W_p$ has simple socle.
\end{prop}

\begin{pf} The proof goes in three steps:
\begin{itemize}
    \item prove that $(V_{p,b_1-p},\cdots, V_{p,b_1})$ is a normal sequence and identify its head;
    \item use Proposition \ref{Linfty}(ii) to prove that $(V(\bomega_{b_1-p, b_1}), V_{p,b_1-p},\cdots, V_{p,b_1})$ is normal;
    \item for $d_r\le k\le b_1-p$, assuming that
    $(V'_{p,k+1},\cdots , V'_{p,b_1-p}, V(\bomega_{b_1-p, b_1}),V_{p,b_1-p},\cdots, V_{p,b_1})$ is normal
    with simple head $H_{k+1}$, prove that $(V'_{p,k}, H_{k+1})$ is unmixed. Proposition
    \ref{Linfty}(ii) then shows that
    $(V'_{p,k},\cdots , V'_{p,b_1-p}, V(\bomega_{b_1-p, b_1}),V_{p,b_1-p},\cdots, V_{p,b_1})$ is normal,
    and taking $k=d_r$ shows that $\bs_p$ is normal, as needed.
\end{itemize}
In Section \ref{pfi} we proved that $(V_{p,b_1-p},\cdots, V_{p,b_1})$ is unmixed and hence also
normal. Since the induction hypothesis applies to $\bos_p$ it follows that
$${\rm hd}(V_{p,b_1-p}\otimes \cdots \otimes  V_{p,b_1})\cong \soc M(\bos_p)\cong V(\bomega_{\bos_p^\dagger}),$$
which completes the first step.
\\\\
For the second step, it suffices by Proposition \ref{Linfty}(ii) to prove that
\begin{equation}\label{ineqSp}
    \lie d (V(\bomega_{b_1-p,b_1})^*, V(\bomega_{\bos_p^\dagger})) = \sum_{k=0}^p \lie d (V(\bomega_{b_1-p,b_1})^*, V_{p,b_1-k}).
\end{equation}
Proposition \ref{kkop}(iii) shows that the left hand side is bounded above by the right hand side,
and so it suffices to prove the reverse inequality. Using the first equality in Proposition
\ref{unmixedpair}(i) and applying Proposition \ref{kkop}(ii) we see that
$$\lie d (V(\bomega_{b_1-p,b_1})^*, V_{p,b_1-k})=0,\ \ 1\le k\le p.$$
Hence we must prove that
$$\lie d (V(\bomega_{b_1-p,b_1})^*, V(\bomega_{\bos_p^\dagger})) \geq \lie d (V(\bomega_{b_1-p,b_1})^*, V_{p,b_1}).$$
Recall from Lemma \ref{mpcount} that we set $N=\sum_{c\in \mathbb Z}\mu_{\bos_p}([c,b_1])$. Let
$c_1\le c_2\le\cdots\le c_N$ be such that $\{[c_1,b_1],\cdots, [c_{N},b_1]\}$ is the multiset of
segments of $\bos_p$ with right end point $b_1$, counted with multiplicity. Write
$$V_{p,b_1} = V(\bomega(c_1,b_1))\otimes\cdots\otimes V(\bomega(c_N, b_1)).$$
Proposition \ref{kkop}(ii) applied to $V_{p,b_1}$ gives
\begin{equation}\label{multm1}\lie d (V(\bomega_{b_1-p,b_1})^*, V_{p,b_1}) = \sum_{j=1}^{N} \lie d (V(\bomega_{b_1-p,b_1})^*,V(\bomega(c_j,b_1)))=N,
\end{equation}
where the last equality follows from Proposition \ref{unmixedpair}(ii).
\\\\
By Lemma \ref{mpcount} we can choose $b_1-p\le \ell_1\leq \ell_2\leq \cdots \leq \ell_{N}$ such that
$\{[\ell_1,b_1],\cdots, [\ell_{N},b_1]\}$ is the multiset of segments of $\bos_p^\dagger$ with right end
point $b_1$, counted with multiplicity. Define $\bomega_1,\bomega_2\in\cal I_n^+$ by
$$\bomega_1 = \prod_{j=1}^{N}\bomega_{\ell_j,b_1}, \ \ \ \bomega_{\bos_p^\dagger} = \bomega_1\bomega_2.$$
By \cite[Corollary 3.2]{BC25b} the module $V(\bomega_1)\otimes V(\bomega_2)$ is $\ell$-highest weight
with head $V(\bomega_{\bos_p^\dagger})$. Moreover, since the pairs $([b_1-p,b_1],[\ell_j,b_1])$ are not
connected, Propositions \ref{chb}(ii) and \ref{kkop}(iv) give
$\lie d (V(\bomega_{b_1-p,b_1}), V(\bomega_{\ell_j,b_1}))=0$ for $1\leq j\leq N$, and hence an
application of Proposition \ref{kkop}(ii) gives
$$ \lie d(V(\bomega_{b_1-p,b_1}),V(\bomega_1))=0.$$
Since $V(\bomega_{b_1-p,b_1})$ is real it follows from \cite[Lemma 4.3(ii)]{KKOP} that
$$\Lambda(V(\bomega_{b_1-p,b_1}), V(\bomega_{\bos_p^\dagger})) = \Lambda(V(\bomega_{b_1-p,b_1}), V(\bomega_1))+\Lambda(V(\bomega_{b_1-p,b_1}), V(\bomega_2)).$$
Using Proposition \ref{Linfty}(i) and \eqref{ladd} we get
$$\lie d (V(\bomega_{b_1-p,b_1})^*, V(\bomega_{\bos_p^\dagger})) = \lie d (V(\bomega_{b_1-p,b_1})^*, V(\bomega_1))+ \lie d (V(\bomega_{b_1-p,b_1})^*, V(\bomega_2)).$$
Since $b_1-p\le\ell_j$ for all $j$, Proposition \ref{kkop}(ii) and Proposition \ref{unmixedpair}(ii) give
\begin{equation*}\lie d (V(\bomega_{b_1-p,b_1})^*, V(\bomega_{\bos_p^\dagger})) \geq \lie d (V(\bomega_{b_1-p,b_1})^*, V(\bomega_1)) = \sum_{j=1}^{N}\lie d (V(\bomega_{b_1-p,b_1})^*, V(\bomega_{\ell_j,b_1})) = N,
\end{equation*}
which together with \eqref{multm1} proves \eqref{ineqSp} and completes the second step.
\\\\
For the third step, let $d_r\leq k\leq b_1-p$ and assume that the sequence
$$(V_{p,k+1}',\cdots, V_{p,b_1-p}', V(\bomega_{b_1-p,b_1}), V_{p,b_1-p},\cdots, V_{p,b_1})$$
is normal; when $k=b_1-p$ this is the sequence $(V(\bomega_{b_1-p,b_1}), V_{p,b_1-p},\cdots, V_{p,b_1})$,
which is normal by the second step. Recall that $H_{k+1}$ is the simple head of the corresponding
tensor product. Since $H_{k+1}$ is a subquotient of this tensor product, Proposition \ref{kkop}(iii)
shows that $(V'_{p,k},H_{k+1})$ is unmixed if $V'_{p,k}$ is unmixed with each of the modules
$$V_{p,s}',\ \ k<s\leq b_1-p, \qquad V_{p,b_1-j},\ \ 0\leq j\leq p, \qquad V(\bomega_{b_1-p,b_1});$$
note that there are no modules of the first kind if $k=b_1-p$. For the first two families this
follows from Propositions \ref{chb}(iii), \ref{kkop}(iv) and \ref{kkop}(i),(ii); note that if
$k=b_1-p=b_1-j$ the factors of $V'_{p,k}$ correspond to intervals $[c,b_1-p]$ with $c<a_p$, while those
of $V_{p,b_1-p}$ correspond to intervals $[c,b_1-p]$ with $c\ge a_p$, so that the second case of
Proposition \ref{chb}(iii) applies. For the last module it follows from Proposition
\ref{unmixedpair}(i), since the factors of $V'_{p,k}$ correspond to intervals with right end point
$k\le b_1-p$. Hence both sides of the equality in Proposition \ref{Linfty}(ii) are zero, which
completes the third step. Finally, since $\bs_p$ is normal, Proposition \ref{simplehead} shows that
$W_p$ has simple socle.
\end{pf}

\subsubsection{} We now prove the final step.
\begin{lem}\label{wm}
We have
$$\bos_m+\bos_m'=\bos[1],\qquad W_m\cong M(\bos[1])\otimes V(\bomega_{b_1-m,b_1}),\qquad
V(\bomega_{\bos^\dagger})\hookrightarrow W_m.$$
In particular, $\soc W_p\cong V(\bomega_{\bos^\dagger})$ for $1\le p\le m$.
\end{lem}
\begin{pf} The equality $\bos_m+\bos_m'=\bos[1]$ is immediate from the definitions.
\\\\
By definition we have $W_m=M(\bos_m)\otimes V(\bomega_{b_1-m,b_1})\otimes M(\bos'_m)$. Recall that
$\mu_{\bos'_m}([c,d])\neq 0$ only if either $d<b_1-m$, or $d=b_1-m$ and $c<a_m$, and the maximality of
$m$ shows that the second case does not occur. Hence Proposition \ref{unmixedpair}(iii) shows that
$V(\bomega_{b_1-m,b_1})$ commutes with $M(\bos'_m)$, and so
$$W_m\cong M(\bos_m)\otimes M(\bos'_m)\otimes V(\bomega_{b_1-m,b_1})\cong M(\bos[1])\otimes
V(\bomega_{b_1-m,b_1}),$$
where the last isomorphism holds since $\bos_m+\bos'_m$ is in lexicographic order.
\\\\
Since $|\bos[1]|<|\bos|$ and $n(\bos[1])\le n(\bos)\le n$, the induction hypothesis gives
$\soc M(\bos[1])\cong V(\bomega_{(\bos[1])^\dagger})$, and hence
$$V(\bomega_{(\bos[1])^\dagger})\otimes V(\bomega_{b_1-m,b_1})\hookrightarrow M(\bos[1])\otimes
V(\bomega_{b_1-m,b_1})\cong W_m.$$
Every interval of $(\bos[1])^\dagger$ has right end point at most $b_1$, and those with right end
point $b_1$ are not connected to $[b_1-m,b_1]$. Hence \cite[Corollary 3.2]{BC25b} shows that
$V(\bomega_{b_1-m,b_1})\otimes V(\bomega_{(\bos[1])^\dagger})$ is $\ell$-highest weight with simple head
$V(\bomega_{b_1-m,b_1}\bomega_{(\bos[1])^\dagger})$, which is isomorphic to $V(\bomega_{\bos^\dagger})$ since
$\bos^\dagger=([b_1-m,b_1])+(\bos[1])^\dagger$ by \eqref{lead}. Since $V(\bomega_{b_1-m,b_1})$ is real,
\cite[Theorem 3.12]{KKKO15} shows that the module on the left of the preceding inclusion has simple
socle isomorphic to this head, i.e.,
$$V(\bomega_{\bos^\dagger}) \cong \soc \big(V(\bomega_{(\bos[1])^\dagger})\otimes
V(\bomega_{b_1-m,b_1})\big)\hookrightarrow W_m,$$
which proves the third statement.
\\\\
Finally, since each $W_p$ has simple socle by Proposition \ref{Spnormal}, the inclusions
\eqref{uw} give $\soc W_p\cong\soc W_m\cong V(\bomega_{\bos^\dagger})$ for $1\le p\le m$.
\end{pf}

\subsection{}\label{socstd} We now prove Proposition \ref{socW}. Let $\bos\in{\rm mult}(\mathbb I_n)$ with $n\ge n(\bos)$. By \cite[Theorem 1]{BC25b},  we have,  $$V(\bomega_{\bos^\sharp})\cong W(\bomega_{\bos^\sharp})\hookrightarrow W(\bomega_\bos)$$ and hence it suffices to show that $W(\bomega_\bos)$ has simple socle. \\\\
Write $\bos=\{[c_1,d_1],\cdots,[c_r,d_r]\}$ with $d_1\ge\cdots\ge d_r$ and suppose first that also  $c_1\geq \cdots\geq c_r$. 
By Proposition \ref{simplehead} it suffices to show that $(V(\bomega_{c_r,d_r}),\cdots, V(\bomega_{c_1,d_1}))$ is normal, and we proceed by induction on $r$, which clearly begins at $r=1$. The induction hypothesis applies to $\bos'= \{[c_1,d_1],\cdots, [c_{r-1},d_{r-1}]\}$  and gives that the sequence
$(V(\bomega_{c_{r-1},d_{r-1}}),\cdots, V(\bomega_{c_1,d_1}))$ is normal and $\soc W(\bomega_{\bos'})\cong V(\bomega_{(\bos')^\sharp})$. By Proposition \ref{Linfty}(ii) the inductive step follows if we show that
\begin{equation}\label{Wnormal}\lie d (V(\bomega_{c_r,d_r})^*, V(\bomega_{(\bos')^\sharp})) = \sum_{j=1}^{r-1} \lie d (V(\bomega_{c_r,d_r})^*, V(\bomega_{c_j,d_j})).
\end{equation}
Recall that $V(\bomega_{c_r,d_r})^* \cong V(\bomega_{d_r,c_r+n+1})$. Since $c_r\leq c_s$ and $d_r\leq d_s<c_r+n+1$ for $1\leq s\leq r$, the pair $([d_r,c_r+n+1],[c_s,d_s])$ is connected if and only if $c_s<d_r$. By Proposition \ref{chb}(ii) together with the fact that $\lie d (V(\bomega_{a,b}),V(\bomega_{c,d}))\leq 1$ for fundamental modules (see \cite{AK97, DO94}), the right hand side of \eqref{Wnormal} is equal to $\#\{j: 1\leq j\leq r-1, \ c_j<d_r\}$.\\\\
On the other hand, by the definition of $(\bos')^\sharp$ there exists $w\in \Sigma_{r-1}$ such that
$$(\bos')^\sharp = ([c_{w(1)},d_1],\cdots, [c_{w(r-1)},d_{r-1}]), \ \ c_{w(s)}\leq d_s, \ \ 1\leq s\leq r-1.$$
 Since $$V(\bomega_{(\bos')^\sharp})\cong \bigotimes_{j=1}^{r-1} V(\bomega_{c_{w(j)},d_j}),$$
by Proposition \ref{kkop}(ii) and the preceding discussion, the left hand side of \eqref{Wnormal} is $\#\{j: c_{w(j)}<d_r\leq d_j\}$, which is equal to the right hand side since $w$ is a permutation. This proves \eqref{Wnormal} and hence Proposition \ref{socW} when $c_1\geq \cdots \geq c_r$. \\\\
For the general case, choose $\sigma\in\Sigma_r$  such that $c_{\sigma(s)}\ge c_{\sigma(s+1)}$, $1\le s<r$, and set $\bos^+=\{[c_{\sigma(1)},d_1],\cdots,[c_{\sigma(r)},d_r]\}\in \mult(\mathbb I_n)$; note that $n(\bos^+)=n(\bos)\leq n$. By \cite[Proposition 1.4.4]{BC25b} we have $(\bos^+)^\sharp =\bos^\sharp$ and, by Proposition  2.3 and Theorem 1 of \cite{BC25b} we have
$$V(\bomega_{\bos^\sharp})\hookrightarrow W(\bomega_\bos)\hookrightarrow W(\bomega_{\bos^+}).$$
The result follows from the previous case, since $\soc W(\bomega_{\bos^+})$ is simple. 

\section{Proof of Theorem \ref{zdualsoc}(ii): the involution $\Phi$}\label{phi}
We first show that $\Phi$ maps the class of a simple module to the class of a simple module, using the generalized quantum affine Schur--Weyl duality functor of Kang--Kashiwara--Kim \cite{KKK15}: on the quiver Hecke side the analogue of $\Phi$ is induced by an algebra involution. The class of the simple module is then identified using the socle statement of Theorem \ref{zdualsoc}(ii), proved in Section \ref{socMs}.

\subsection{}\label{window} For a fixed $\bos = \{[c_1,d_1],\cdots, [c_r,d_r]\}\in{\rm mult}(\mathbb I_n)$
satisfying $n\ge n(\bos)$, let $c=\min_j c_j$ and write $\bou=\bou_c$ (see Section \ref{bou}). Since
$c\leq c_k<d_k\le c+n$ for all $1\leq k\leq r$, the module $W(\bomega_\bos)$ is an object of
$\mathscr F_n(\bou)$, and hence so is every composition factor of $W(\bomega_\bos)$. Since
$\Phi([V(\bomega_{a,b})])=[V(\bomega(a,b))]$ and $V(\bomega(a,b))$ is an object of $\mathscr F_n(\bou)$ if
$c\le a<b\le c+n$, it follows that $\Phi$ restricts to an involution of $\cal K(\mathscr F_n(\bou))$.
\\\\
The category $\mathscr F_n(\bou)$ contains exactly $n(n+1)/2$ fundamental modules, namely the modules
$V(\bomega_{a,b})$ with $c\le a<b\le c+n$, and its simple objects are the simple constituents of tensor
products of these modules. For each $1\le i\le n$ the fundamental $\ell$-weights $\bomega_{a,a+i}$,
$c\le a\le c+n-i$, correspond to the $\ell$-weights $Y_{i,q^{2c+i+2k}}$, $0\le k\le n-i$, in the
notation of \cite{FR99}. Comparing with \cite[Example 3.8, (3.3) and Remark 3.6]{HL15}, we see that,
up to a shift of the spectral parameter, these are exactly the highest $\ell$-weights of the
fundamental modules of the category $\mathscr C_Q$ attached to the equioriented quiver $Q$ of type
$A_n$. Since both categories are the full subcategories whose simple constituents are the simple
constituents of tensor products of their fundamental modules, it follows that
$\mathscr F_n(\bou)=\mathscr C_Q$.

\subsection{Quiver Hecke algebras}\label{klr} Recall that for $\beta=\sum_{i}m_i\alpha_i$ of height $m=\sum_im_i$, the algebra $R(\beta)$ is generated
by idempotents $e(\nu)$, $\nu=(\nu_1,\cdots,\nu_m)$ a sequence of vertices with
$\alpha_{\nu_1}+\cdots+\alpha_{\nu_m}=\beta$, by elements $x_1,\cdots,x_m$ and by elements
$\tau_1,\cdots,\tau_{m-1}$, subject to the relations of \cite[Section 1.2]{KKK18} (see also
\cite{KL09}). Let $R=\bigoplus_\beta R(\beta)$ be the quiver Hecke algebra of type $A_n$ with vertices $c,c+1,\cdots,c+n-1$ and the parameters of \cite[(2.4)]{KKK15}, and let $R$-gmod be the category of finite-dimensional graded $R$-modules, with the convolution product $\circ$. Here and in Section \ref{sec:sigma}, $q$ also denotes the grading shift $M\mapsto qM$ on $R$-gmod.\\\\
For $c\le a<b\le c+n$ write $\alpha_{a,b}=\alpha_a+\cdots+\alpha_{b-1}$, and let $L^+_{a,b}$ (resp. $L^-_{a,b}$) be the one-dimensional $R(\alpha_{a,b})$-module on which the $x_k$ and $\tau_k$ act by zero and $e(\nu)$ acts by $1$ if $\nu=(a,a+1,\cdots,b-1)$ (resp. $\nu=(b-1,\cdots,a+1,a)$) and by $0$ otherwise, as in \cite[(4.2.1)]{KKK18}. 
\\\\
In \cite{KKK15} the authors defined a functor $\mathcal F:R\text{-gmod}\to\mathscr F_n(\bou)$ which has the following properties.
\begin{prop}\label{KLR}
\begin{enumerit}
\item[(i)] $\mathcal F$ is exact, $\mathcal F(M\circ N)\cong\mathcal F(M)\otimes\mathcal F(N)$, and $\mathcal F(qM)\cong\mathcal F(M)$ \cite[Theorem 2.1.1]{KKK15}.
\item[(ii)] $\mathcal F(L_{j,j+1}^+)\cong V(\bomega_{j,j+1})$ for $c\le j<c+n$ \cite[(4.7)]{KKK15}.
\item[(iii)] $\mathcal F$ induces an isomorphism $\cal K(R\text{-gmod})/(q-1)\cong\cal K(\mathscr F_n(\bou))$ which maps the classes of simple modules bijectively onto the classes of simple modules \cite[Theorem 4.2.1, Corollaries 4.3.5 and 4.3.6]{KKK15}.
\end{enumerit}
\end{prop}
\subsubsection{}
For $\beta$ of height $m$ let $\sigma$ be the automorphism of $R(\beta)$ given by
$$e(\nu)\mapsto e(\nu^{\rm rev}),\qquad x_k\mapsto x_{m+1-k},\qquad \tau_ke(\nu)\mapsto(-1)^{\delta_{\nu_k,\nu_{k+1}}}\tau_{m-k}e(\nu^{\rm rev}),$$
where $\nu^{\rm rev}=(\nu_m,\cdots,\nu_1)$; see \cite[Section 2.1]{KL09} and \cite[Section 1.1.4, (1.32)]{LV11}. Twisting the action by $\sigma$ defines an exact autoequivalence $M\mapsto M^\sigma$ of $R$-gmod which maps simple modules to simple modules  and induces an isomorphism $(M\circ N)^\sigma\cong N^\sigma\circ M^\sigma$. Since $\cal K(\mathscr F_n(\bou))$ is commutative, Proposition \ref{KLR}(iii) shows that $\sigma$ induces a ring involution $\Omega$ of $\cal K(R\text{-gmod})/(q-1)$ with $\Omega([M])=[M^\sigma]$, which maps classes of simple modules to classes of simple modules. Clearly $(L^+_{a,b})^\sigma\cong L^-_{a,b}$.
\begin{lem}\label{Lpm}
For $c\le a<b\le c+n$ we have $\mathcal F(L^+_{a,b})\cong V(\bomega_{a,b})$ and $\mathcal F(L^-_{a,b})\cong V(\bomega(a,b))$.
\end{lem}
\begin{pf}
By Frobenius reciprocity there are surjections, up to grading shifts,
$$L_{a,a+1}^+\circ L_{a+1,a+2}^+\circ\cdots\circ L_{b-1,b}^+\to L^+_{a,b},\qquad L_{b-1,b}^+\circ\cdots\circ L_{a+1,a+2}^+\circ L_{a,a+1}^+\to L^-_{a,b}.$$
By Proposition \ref{KLR}, $\mathcal F(L^+_{a,b})$ and $\mathcal F(L^-_{a,b})$ are simple quotients of
$$V(\bomega_{a,a+1})\otimes\cdots\otimes V(\bomega_{b-1,b})\qquad\text{and}\qquad V(\bomega_{b-1,b})\otimes\cdots\otimes V(\bomega_{a,a+1})$$
respectively. The second module is $M(\boT)$, with $\boT=\{[b-1,b],\cdots, [a,a+1]\}$. By Theorem \ref{zdualsoc} its head is $V(\bomega(a,b))$ and its socle is $V(\bomega_{a,b})$, since $\boT^\dagger = \{[a,b]\}$. Hence, by Proposition \ref{kkop}(v), the head of the first module is $V(\bomega_{a,b})$.
\end{pf}
\subsection{} Let $\varpi$ denote the isomorphism of rings induced by $\mathcal F$ in Proposition \ref{KLR}(iii).
\begin{prop}\label{Phistar}
We have $\Phi\circ\varpi=\varpi\circ\Omega$ on $\cal K(R\text{-gmod})/(q-1)$. In particular, $\Phi$ maps the class of every simple object of $\mathscr F_n(\bou)$ to the class of a simple object of $\mathscr F_n(\bou)$.
\end{prop}
\begin{pf}
Both sides are ring homomorphisms. By Lemma \ref{Lpm}, we have
$$\Phi\varpi([L^+_{a,b}])=\Phi([V(\bomega_{a,b})])=[V(\bomega(a,b))]=\varpi([L^-_{a,b}])=\varpi\Omega([L^+_{a,b}]),$$ and the classes $[L^+_{a,b}]$ generate $\cal K(R\text{-gmod})/(q-1)$, since by Proposition \ref{KLR}(iii) and Lemma \ref{Lpm} they correspond to the generators $[V(\bomega_{a,b})]$ of $\cal K(\mathscr F_n(\bou))$. The second assertion follows from Proposition \ref{KLR}(iii).
\end{pf}

\subsection{Completion of the proof of Theorem \ref{zdualsoc}(ii)}
Recall the expansion \eqref{Wexp}.
\begin{prop}\label{OmegaPhi}
Let $\bos\in{\rm mult}(\mathbb I_n)$ with $n\ge n(\bos)$. Then $\Phi([V(\bomega_\bos)])=[V(\bomega_{\bos^\dagger})]$.
\end{prop}
\begin{pf}
We proceed by induction on $\lambda_\bos$ with respect to the dominance order. By Section \ref{window} and the induction hypothesis, $\Phi([V(\bomega_{\mathbf t})])=[V(\bomega_{\mathbf t^\dagger})]$ whenever $m_{\bos,\mathbf t}\neq0$. Since $\Phi([W(\bomega_\bos)])=[M(\bos)]$, applying $\Phi$ to \eqref{Wexp} gives
$$[M(\bos)]=\Phi([V(\bomega_\bos)])+\sum_{\mathbf t}m_{\bos,\mathbf t}\,[V(\bomega_{\mathbf t^\dagger})].$$
By the socle statement of Theorem \ref{zdualsoc}(ii), $V(\bomega_{\bos^\dagger})$ is a composition factor of $M(\bos)$. Since the map $\bos\mapsto\bos^\dagger$ is injective and $\mathbf t\neq\bos$, the module $V(\bomega_{\bos^\dagger})$ is not among the $V(\bomega_{\mathbf t^\dagger})$. Since all coefficients are non-negative and $\Phi([V(\bomega_\bos)])$ is the class of a simple module by Proposition \ref{Phistar}, we get $\Phi([V(\bomega_\bos)])=[V(\bomega_{\bos^\dagger})]$, which completes the proof of Theorem \ref{zdualsoc}(ii).
\end{pf}
\begin{rem}
By Propositions \ref{Phistar} and \ref{OmegaPhi}, if $L$ is a simple $R$-module with $\mathcal F(L)\cong V(\bomega_\bos)$, then $\mathcal F(L^\sigma)\cong V(\bomega_{\bos^\dagger})$: on the quiver Hecke side, $\sigma$ realizes the Zelevinsky involution.
\end{rem}

\section{The autoequivalence $\Sigma$ and $t$-analogues}\label{sec:sigma}
In this section we explain how the involution $\Phi$ is induced, on each of the categories
$\mathscr F_n(\bou_c)$, by an exact autoequivalence, and how the identity of Section \ref{multdual} lifts to
quantum Grothendieck rings.

\subsection{The autoequivalence $\Sigma$}\label{Sigma} Let $c\in\mathbb Z$ and write $\bou=\bou_c$ (see Section
\ref{bou}). The functor $\mathcal F$ of Section \ref{klr} is defined on the category $R$-gmod of finite-dimensional
graded $R$-modules, and since $\mathcal F(qM)\cong\mathcal F(M)$ it only induces an isomorphism of Grothendieck rings
after passing to $q=1$ (Proposition \ref{KLR}(iii)). Let $R$-fmod be the category of finite-dimensional, not
necessarily graded, $R$-modules on which the center of $R$ acts nilpotently, and let $\mathcal F^u:R\text{-fmod}\to\mathscr F_n(\bou)$
be the functor defined in the same way as $\mathcal F$, so that $\mathcal F\cong\mathcal F^u\circ{\rm Forget}$, where
${\rm Forget}:R\text{-gmod}\to R\text{-fmod}$ forgets the grading. It was proved in \cite[Theorem B]{Fuj} that $\mathcal F^u$ is an
equivalence of monoidal categories.
\\\\
Since every simple object of $R$-fmod is gradable, uniquely up to a grading shift \cite[Section 2.5]{KL09}, ${\rm Forget}$ induces
an isomorphism $\cal K(R\text{-gmod})/(q-1)\cong\cal K(R\text{-fmod})$ which maps classes of simple modules bijectively onto
classes of simple modules. The automorphism $\sigma$ of Section \ref{klr} is a graded automorphism, so twisting by
$\sigma$ is an exact autoequivalence of $R$-fmod which commutes with ${\rm Forget}$, and we obtain an exact
autoequivalence
$$\Sigma=\mathcal F^u\circ(-)^\sigma\circ(\mathcal F^u)^{-1}:\mathscr F_n(\bou)\to\mathscr F_n(\bou).$$
Under the isomorphisms above, the map induced by $\Sigma$ on $\cal K(\mathscr F_n(\bou))$ corresponds to the involution
$\Omega$ of $\cal K(R\text{-gmod})/(q-1)$, and hence by Proposition \ref{Phistar} it is $\Phi$.
\begin{prop}\label{Sigmaprop} The autoequivalence $\Sigma$ of $\mathscr F_n(\bou_c)$ has the following properties.
\begin{enumerit}
\item[(i)] For all objects $M,N$ of $\mathscr F_n(\bou_c)$ we have $\Sigma(M\otimes N)\cong\Sigma(N)\otimes\Sigma(M)$,
and $\Sigma^2$ is isomorphic to the identity functor.
\item[(ii)] For every simple object $V(\bomega_\bos)$ of $\mathscr F_n(\bou_c)$ we have
$\Sigma(V(\bomega_\bos))\cong V(\bomega_{\bos^\dagger})$. In particular
$\Sigma(V(\bomega_{a,b}))\cong V(\bomega(a,b))$ for $c\le a<b\le c+n$.
\item[(iii)] For all objects $M$ of $\mathscr F_n(\bou_c)$ we have $[\Sigma(M)]=\Phi([M])$, and
$$\soc\Sigma(M)\cong\Sigma(\soc M),\qquad {\rm hd}\,\Sigma(M)\cong\Sigma({\rm hd}\,M).$$
\end{enumerit}
\end{prop}
\begin{proof}[Proof of Proposition \ref{Sigmaprop}] Part (i) follows since $\mathcal F^u$ is an equivalence of monoidal
categories, twisting by $\sigma$ satisfies $(M\circ N)^\sigma\cong N^\sigma\circ M^\sigma$ in $R$-fmod, and $\sigma^2={\rm id}$. The
first statement of (iii) is the last sentence of the preceding discussion, and the second holds since $\Sigma$ is
an exact autoequivalence. For (ii), $\Sigma$ maps simple objects to simple objects, so by (iii)
$[\Sigma(V(\bomega_\bos))]=\Phi([V(\bomega_\bos)])$ is the class of a simple object; since every simple object of
$\mathscr F_n(\bou)$ satisfies $n\ge n(\bos)$, Proposition \ref{OmegaPhi} gives $\Sigma(V(\bomega_\bos))\cong V(\bomega_{\bos^\dagger})$.
The second statement of (ii) is the case $\bos=([a,b])$, or directly Lemma \ref{Lpm}.
\end{proof}
\noindent Notice that $\Sigma(W(\bomega_\bos))$ is the tensor product of the modules $V(\bomega(c_s,d_s))$ in
the order opposite to that of $M(\bos)$, so that Proposition \ref{Sigmaprop}(iii) and Theorem
\ref{zdualsoc} are compatible: the head of $\Sigma(W(\bomega_\bos))$ is $\Sigma(V(\bomega_\bos))\cong
V(\bomega_{\bos^\dagger})$, and by Proposition \ref{kkop}(v) the head of a tensor product of the real modules
$V(\bomega(c_s,d_s))$ in the opposite order is the socle of $M(\bos)$.

\subsection{$t$-analogues}\label{tsec} In this subsection we show that the identity of Section \ref{multdual} admits a
natural $t$-deformation, in which Grothendieck rings are replaced by quantum Grothendieck rings. The idea is the
following. The involution $\Omega$ of Section \ref{klr} is induced by the automorphism $\sigma$ of $R$, and since $\sigma$
preserves the grading it makes sense before setting $q=1$. On the other hand, by the work of Hernandez--Leclerc
\cite{HL15} and Kang--Kashiwara--Kim \cite{KKK15}, the graded Grothendieck ring of $R$ is, up to an explicit twist of the
product, the quantum Grothendieck ring of $\mathscr F_n(\bou_c)$. Transporting the twist by $\sigma$ along this
identification gives the required $t$-analogue of $\Phi$.\\\\
 Fix $c\in\mathbb Z$ and write $\bou=\bou_c$. Let $\cal K_t(\mathscr F_n(\bou))$ be the
$\mathbb Z[t^{\pm1/2}]$-span of the $(q,t)$-characters $[V]_t$ of the simple objects $V$ of $\mathscr F_n(\bou)$ in the quantum
Grothendieck ring \cite{Nak04,VV03,Her04}, normalized as in \cite[Section 5]{HL15}, and denote its product by $*$. Since
$\mathscr F_n(\bou)$ is a shift of $\mathscr C_Q$ (Section \ref{window}), this is a $\mathbb Z[t^{\pm\frac12}]$-form of the ring $\cal K_{t,Q}$ of \cite[Section 5.12]{HL15}.
Its structure constants in the basis $\{[V]_t\}$ lie in $\mathbb N[t^{\pm1/2}]$ \cite{VV03}, and $t^{1/2}\mapsto1$ maps $[V]_t$ to
$[V]$ \cite[Theorem 5.2]{HL15}.
\\\\
Let $\psi$ be the anti-involution of $R$ fixing $e(\nu),x_k,\tau_k$, which defines the duality $M\mapsto M^\circledast$ on $R$-gmod.
Since $\sigma$ preserves degrees and commutes with $\psi$, twisting by $\sigma$ induces a $\mathbb Z[q^{\pm1}]$-linear anti-involution
$\Omega_q$ of the graded Grothendieck ring $\cal K(R\text{-gmod})$, which permutes the classes of the self-dual simple modules and
specializes to $\Omega$ at $q=1$. By \cite[Theorem 6.1]{HL15}, \cite[Corollary 4.3.5]{KKK15} and \cite{KL09,VV11}, there is a bijective
$\mathbb Z[t^{\pm1/2}]$-linear map
$$\varpi_t:\cal K(R\text{-gmod})\otimes_{\mathbb Z[q^{\pm1}]}\mathbb Z[t^{\pm1/2}]\to\cal K_t(\mathscr F_n(\bou)),\qquad q\mapsto t^{-1},$$
with the following properties. It maps the class of every self-dual simple module $L$ to $[\mathcal F(L)]_t$. For $x,y$ of
weights $\beta,\gamma$ it satisfies $\varpi_t(xy)=t^{(\beta,\gamma)/2}\varpi_t(x)*\varpi_t(y)$, where $(\cdot,\cdot)$ is the
symmetric form on the root lattice of type $A_n$. (The factor comes from the rescaling by $t^{N(\beta)/2}$,
$N(\beta)=(\beta,\beta)/2-{\rm ht}\,\beta$, in \cite[Theorem 6.1]{HL15}.) Set $\Phi_t=\varpi_t\circ\Omega_q\circ\varpi_t^{-1}$.

\begin{prop}\label{Phit} $\Phi_t$ is a $\mathbb Z[t^{\pm1/2}]$-linear anti-involution of $\cal K_t(\mathscr F_n(\bou))$ which
specializes to $\Phi$ at $t=1$, and $\Phi_t([V(\bomega_\bos)]_t)=[V(\bomega_{\bos^\dagger})]_t$ for every simple object
$V(\bomega_\bos)$ of $\mathscr F_n(\bou)$.
\end{prop}
\begin{pf} The specialization is $\varpi\Omega\varpi^{-1}=\Phi$ by Proposition \ref{Phistar}. Since $\Omega_q$ reverses products
and preserves weights, and $t^{(\beta,\gamma)/2}$ is symmetric in $\beta,\gamma$, $\Phi_t$ reverses products.\\\\
Let $L$ be a
self-dual simple module with $\mathcal F(L)\cong V(\bomega_\bos)$ (see \cite{KL09}). Then $L^\sigma$ is simple and self-dual
 and, since $n\ge n(\bos)$,
$\mathcal F(L^\sigma)\cong V(\bomega_{\bos^\dagger})$ by the Remark after Proposition \ref{OmegaPhi}. Hence
$\Phi_t([V(\bomega_\bos)]_t)=\varpi_t([L^\sigma])=[V(\bomega_{\bos^\dagger})]_t$.
\end{pf}
\begin{cor}\label{multdualt} Let $\bos=\{[c_1,d_1],\cdots,[c_r,d_r]\}$ with $d_1\ge\cdots\ge d_r$ and $n\ge n(\bos)$, and take
$c=\min_jc_j$. Define $P_{\bos,\boT}(t)$ by
$$[V(\bomega_{c_1,d_1})]_t*\cdots*[V(\bomega_{c_r,d_r})]_t=\sum_\boT P_{\bos,\boT}(t)\,[V(\bomega_\boT)]_t.$$
Then $P_{\bos,\boT}(t)\in\mathbb N[t^{\pm1/2}]$, $P_{\bos,\boT}(1)=m_{\bos,\boT}$, and
$$[V(\bomega(c_r,d_r))]_t*\cdots*[V(\bomega(c_1,d_1))]_t=\sum_\boT P_{\bos,\boT}(t)\,[V(\bomega_{\boT^\dagger})]_t.$$
\end{cor}
\begin{pf} Positivity follows from that of the structure constants, and $P_{\bos,\boT}(1)=m_{\bos,\boT}$ follows from
\eqref{Wexp} by specializing at $t=1$. The last identity follows by applying $\Phi_t$, using Proposition \ref{Phit} and
$\{[a,b]\}^\dagger=\{[a,a+1],\cdots,[b-1,b]\}$.
\end{pf}
\noindent At $t=1$ this is the identity of Section \ref{multdual}. On the quiver Hecke side, the twist by $\sigma$ is Kashiwara's
involution on the dual canonical basis \cite{KL09,LV11}; for its relation with the Zelevinsky involution see \cite{LTV99,JL10}.

\subsection{Consequences for quiver Hecke algebras}\label{klrcons}
Retain the notation of Sections \ref{klr} and \ref{Sigma}.\\\\ The positive roots of $R$ are  $\alpha_{a,b}$, $c\le a<b\le c+n$. A  multiset of positive roots corresponds in the obvious way to a multisegment $\bos$ all of whose segments are contained in $[c,c+n]$, in particular we have $n\ge n(\bos)$. For $\bos=\{[c_1,d_1],\cdots,[c_r,d_r]\}$, written in the order of Section \ref{ms}, set
\begin{gather*}
\Delta^+(\bos)=L^+_{c_1,d_1}\circ\cdots\circ L^+_{c_r,d_r},\qquad \nabla^-(\bos)=L^-_{c_1,d_1}\circ\cdots\circ L^-_{c_r,d_r},\\
\nabla^+(\bos)=\nabla^-(\bos)^\sigma\cong L^+_{c_r,d_r}\circ\cdots\circ L^+_{c_1,d_1}\qquad \Delta^-(\bos)=\Delta^+(\bos)^\sigma\cong L^-_{c_r,d_r}\circ\cdots\circ L^-_{c_1,d_1}.
\end{gather*}
By Proposition \ref{KLR}(i), Lemma \ref{Lpm} and the definition of $\Sigma$,
\begin{gather}\label{Fstd}
\mathcal F(\Delta^+(\bos))\cong W(\bomega_\bos),\quad \mathcal F(\nabla^-(\bos))\cong M(\bos),\\
% \label{Fcostd} \mathcal F(\nabla^+(\bos))\cong \Sigma(M(\bos))\cong M(\bos)^{\op}\qquad \cal F(\Delta^-(\bos))\cong \Sigma(W(\bomega_\bos))\cong W(\bomega_\bos)^{\op},\\
\label{Fcostd} \mathcal F(\nabla^+(\bos))\cong \Sigma(M(\bos))\cong W(\bomega_{\bos})^{\op},\qquad \cal F(\Delta^-(\bos))\cong \Sigma(W(\bomega_\bos))\cong M(\bos)^{\op},
\end{gather}
where $(V_1\otimes \cdots\otimes V_k)^{\op}\cong V_k\otimes\cdots\otimes  V_1$.\\\\
Using Proposition \ref{KLR}(iii) and Proposition \ref{Sigmaprop}(ii) we see that there exist
 $L^\pm(\bos)$  simple graded $R$-modules (unique up to grading shift), with $$\mathcal F(L^+(\bos))\cong V(\bomega_\bos),\ \ L^-(\bos)=L^+(\bos)^\sigma,\ \ \mathcal F(L^-(\bos))\cong V(\bomega_{\bos^\dagger}).$$ 

\begin{prop}\label{klrsoc} Let $\bos$ be a multisegment all of whose segments are contained in $[c,c+n]$. Up to grading shifts:
\begin{enumerit}
\item[(i)] $\Delta^\pm(\bos)$ has simple head $L^\pm(\bos)$ and simple socle $L^\pm(\bos^\sharp)$;
\item[(ii)] $\nabla^\pm(\bos)$ has simple socle $L^\pm(\bos)$ and simple head $L^\pm(\bos^\sharp)$.
\end{enumerit}
Moreover, \begin{gather*} L^+(\bos^\sharp) \cong \Delta^+(\bos^\sharp),\ \ L^-(\bos^\sharp)\cong\nabla^-(\bos^\sharp),\ \
\mathcal F(L^-(\bos^\sharp))\cong V(\bomega(\bos)).
\end{gather*}

\end{prop}

\begin{pf} Recall that $\mathcal F\cong\mathcal F^u\circ{\rm Forget}$. Since  $\mathcal F^u$ is an equivalence the  socle and head  of a graded modules is graded. Moreover, a graded module is simple if and only if it is simple as an ungraded module \cite[Section 2.5]{KL09}.  It follows that $\mathcal F$ preserves socles and heads and also  that two graded simple modules with isomorphic images under $\mathcal F$ agree up to grading shift. \\\\
Using \eqref{Fstd} and the well-known  fact that $W(\bomega_\bos)$
has simple head shows that $\Delta^+(\bos)$ has simple head. Proposition \ref{socW} shows that $W(\bomega_{\bos})$ has simple socle isomorphic to $W(\bomega_{\bos^\sharp})$ and hence we also have that $\Delta^+(\bos)$ has simple socle $L^+(\bos^\sharp)$. Since the head and socle of $W(\bomega_\bos)$ and $W(\bomega_{\bos})^{\op}$ are interchanged part (ii) follows for $\nabla^+(\bos)$.
\\\\
Similarly, using Theorem \ref{zdualsoc} and  we see that socle of $\nabla^-(\bos)$ is $L^-(\bos)$ and simple head $L^-(\bos^\sharp)$, and
Theorem \ref{zdualsoc}(ii) gives
$$\mathcal F(L^-(\bos^\sharp))\cong V(\bomega_{(\bos^\sharp)^\dagger})\cong V(\bomega(\bos)).$$
Since the head and socle of $M(\bos)$ and $M(\bos)^{\op}$ are interchanged, \eqref{Fcostd} shows that part (i) follows for
$\Delta^-(\bos)$.
\\\\
Finally, $\mathcal F(\Delta^+(\bos^\sharp))\cong W(\bomega_{\bos^\sharp})\cong V(\bomega_{\bos^\sharp})$ and
$\mathcal F(\nabla^-(\bos^\sharp))\cong M(\bos^\sharp)\cong V(\bomega(\bos))$ are simple by Section \ref{std} and
Theorem \ref{zdualsoc}(i), and hence $\Delta^+(\bos^\sharp)\cong L^+(\bos^\sharp)$ and $\nabla^-(\bos^\sharp)\cong L^-(\bos^\sharp)$.

\end{pf}

\begin{rem}  The simple head of $\Delta^\pm(\bos)$ and simple socle of $\nabla^\pm(\bos)$ are part of the general theory of proper standard modules; the simple socle of $\Delta^\pm(\bos)$ and head of $\nabla^\pm(\bos)$ are not, and correspond to Proposition \ref{socW} and Theorem \ref{zdualsoc}(i). In addition, Theorem \ref{zdualsoc}(ii) adds  the identification $\mathcal F(L^-(\bos))\cong V(\bomega_{\bos^\dagger})$, i.e., twisting by $\sigma$ acts on labels of simple modules by the Zelevinsky involution (compare \cite{JL10,LTV99}). Choosing $n$ large, the proposition also applies in type $A_\infty$.

\end{rem}

\end{document}